\documentclass[12pt]{amsart}
\usepackage{hyperref}
\usepackage{amsfonts}
\usepackage{amsmath}
\usepackage{amssymb}
\usepackage{amscd}
\usepackage{graphics}
\usepackage{latexsym}
\usepackage{color}
\usepackage{epsfig}
\usepackage{tikz-cd}
\usepackage{amsopn}

\newtheorem{theorem}{Theorem}[section]
\newtheorem{proposition}[theorem]{Proposition}
\newtheorem{lemma}[theorem]{Lemma}
\newtheorem{corollary}[theorem]{Corollary}

\newtheorem{problem}[theorem]{Problem}

\newtheorem{claim}[theorem]{Claim}

\theoremstyle{remark}

\newtheorem{notation}[theorem]{Notation}

\newtheorem{example}[theorem]{Example}

\newtheorem{definition}[theorem]{Definition}

\newcommand{\PP}{\mathbb{P}}

\newcommand{\CC}{\mathbb{C}}

\newcommand{\ZZ}{\mathbb{Z}}

\newcommand{\Ker}{\mbox{Ker}}

\DeclareMathOperator{\mw}{mw}

\DeclareMathOperator{\PGL}{\mathbb{P}GL}

\DeclareMathOperator{\Stab}{Stab}

\begin{document}

\title{$\PGL$ orbits in tree varieties}
\author[I. Coskun]{Izzet Coskun}
\address{Department of Mathematics, Stat. and CS \\University of Illinois at Chicago, Chicago, IL 60607}
\email{icoskun@uic.edu}
\author[D. Eken]{Demir Eken}
\address{Department of Mathematics, University of Michigan, Ann Arbor, MI 48109}
\email{eken@umich.edu}

\author[C. Yun]{Chris Yun}
\address{Department of Mathematics, Stat. and CS \\University of Illinois at Chicago, Chicago, IL 60607}
\email{cyun2@uic.edu}

\thanks{During the preparation of this article the first author was partially supported by the NSF FRG grant DMS 1664296 and NSF grant DMS-2200684..}
\begin{abstract}
In this paper, we introduce tree varieties as a natural generalization of products of partial flag varieties. We study orbits of the $\PGL$ action on tree varieties. We  characterize tree varieties with finitely many $\PGL$ orbits, generalizing a celebrated theorem of Magyar, Weyman and Zelevinsky. We give criteria that guarantee that a tree variety has a dense $\PGL$-orbit and provide many examples of tree varieties that do not have dense $\PGL$ orbits. We show that a triple of two-step flag varieties $F(k_1, k_2; n)^3$ has a dense $\PGL(n)$ orbit if and only if $k_1 + k_2 \not= n$.
\end{abstract}
\maketitle

\section{Introduction}

In this paper, we introduce tree varieties and study $\PGL$ orbits in tree varieties.  These varieties arise naturally when studying $\PGL$ orbits on products of flag varieties via inductive constructions. We characterize tree varieties with finitely many $\PGL$ orbits, generalizing a celebrated theorem of Magyar, Weyman and Zelevinsky \cite{MWZ}. 

We also study  tree varieties with dense $\PGL$ orbits.  The study of the product of flag varieties with dense $\PGL$ orbits was initiated by Popov \cite{Popov, Popov2} based on a question of M. Burger and further studied in \cite{CHZ, Devyatov}. We refer the reader to \cite{Smirnov} for a recent survey. We give some criteria that guarantee that the tree varieties have dense $\PGL$ orbits. Unfortunately, at present a complete characterization of tree varieties with dense $\PGL$ orbits seems our of reach, even in the special case of  products of three flag varieties. We do, however, settle the first non-trivial case by showing that $F(k_1, k_2; n)^3$ has a dense $\PGL(n)$ orbit if and only if $k_1 + k_2 \not= n$. 

We now introduce tree varieties. Throughout the paper we work over an algebraically closed field of arbitrary characteristic.

\subsection{Tree varieties} Let $T$ be a directed tree. We will denote the vertices of $T$ by $V(T)$ and the edges of $T$ by $E(T)$. Each directed edge is determined by specifying a pair of vertices $(s,t)$, where the edge points from  the {\em source} $s$ to the {\em target} $t$.  

\begin{definition}[Labeled tree]
A {\em labeled tree} $(T, \phi)$ is a pair such that 
\begin{itemize}
\item $T$ is a rooted, directed tree where all the edges point towards the root; and 
\item $\phi: V(T) \to \ZZ_{>0}$ is a function that assigns to each vertex of $T$ a positive integer such that if $(s,t) \in E(T)$, then $\phi(s) < \phi(t)$.
\end{itemize}
\end{definition}
Let $r$ denote the root of the tree $T$. The vertex $r$ is the only vertex which is not the source of an edge. We say $n= \phi(r)$ is the {\em ambient dimension} of $(T, \phi)$. A {\em leaf} of $T$ is a vertex $s$ which is not the target of any edge. A {\em branch} of a tree is a maximal directed chain containing a leaf such that each vertex in the chain is the target of at most one edge. Note that each branch contains a unique leaf. We will depict labeled trees by drawing a tree where each vertex $v$ is labeled by $\phi(v)$. 

\begin{example} \label{ex-tree}
In the tree below, the leaves are the vertices marked by $d_1$ and $d_5$. The tree has two branches, one with three vertices labeled $d_1, d_2, d_3$ and one with one vertex labeled $d_5$. 
\[
\begin{tikzcd}
d_1 \arrow[r] & d_2 \arrow[r] & d_3 \arrow[r]& d_4 \arrow[r] & n \\
& &  & d_5 \arrow[u] & 
\end{tikzcd}
\]

\end{example}

\begin{definition}[Tree variety]\label{def-treevariety}
Let  $(T, \phi)$ be a labeled tree with ambient dimension $n$.  Let $W$ be an $n$-dimensional vector space. The {\em tree variety} $F(T, \phi)$ associated to $(T, \phi)$  is the variety which parameterizes a   $\phi(v)$-dimensional subspace $U_v$ of $W$ for each vertex $v \in V(T)$  such that  $U_s \subset U_t$ whenever $(s,t) \in E(T)$.
\end{definition}

Tree varieties are smooth, irreducible, projective varieties and their dimensions are readily computed (see Theorem \ref{thm-irredanddim}).

\begin{example}\label{ex-flag}
Let  $(T, \phi)$ be the labeled tree where $T$ is a chain with $m+1$ vertices and $\phi$ associates the positive integers $k_1 < k_2 < \cdots < k_m <n$ to these vertices
\[
\begin{tikzcd}
k_1 \arrow[r] & k_2 \arrow[r] & \cdots  \arrow[r] & k_m \arrow[r] & n
\end{tikzcd}
\]
In this case, the tree variety $F(T, \phi)$ is the $m$-step partial flag variety $F(k_1, \dots, k_m; n)$ parameterizing partial flags $U_1 \subset U_2 \subset \cdots \subset U_m \subset W$, where $U_i$ has dimension $k_i$. 
\end{example}

\begin{example}\label{ex-productofflags}
Consider the labeled tree  $(T, \phi)$, where $T$ is a union of $j$ chains joined at the root.  
\[
\begin{tikzcd}
k_{1,1} \arrow[r] & k_{1,2} \arrow[r] & \cdots  \arrow[r] & k_{1, m_1} \arrow[rd] & \\ 
& & \cdots & &  n \\
k_{j,1} \arrow[r] & k_{j,2} \arrow[r] & \cdots  \arrow[r] & k_{j, m_j} \arrow[ru] &  
\end{tikzcd}
\]
In this case, the tree variety $F(T, \phi) = \prod_{i=1}^j F(k_{i,1}, \dots, k_{i,m_i}; n)$, the product of $j$-partial flag varieties. Hence, tree varieties generalize products of partial flag varieties. They are also closely related to quiver varieties associated to a rooted, directed tree. However, unlike in quiver varieties, in tree varieties we do not take any quotients. 
\end{example}

If the ambient dimension of $(T, \phi)$ is $n$, the group $\PGL(n)$ acts on $F(T, \phi)$. In this paper, we are interested in the orbits of this action. We address the following two main questions.

\begin{enumerate}
\item When does the action of $\PGL(n)$  on $F(T, \phi)$ have finitely many orbits?
\item When does the action of $\PGL(n)$  on $F(T, \phi)$ have a dense orbit?
\end{enumerate}
We resolve the first of these questions completely. The second question is much harder, nevertheless, we obtain many new partial results. In fact, the main motivation for introducing tree varieties came from studying the second question for products of Grassmannians.

\subsection{Results} We now describe our results in detail.

\begin{definition}
Given a branch $B$ of a labeled tree $(T, \phi)$, let $s_B$ denote the leaf of $B$. Then the  {\em minimum width} $\mw(B)$ of $B$ is defined by $$\mw(B) := \min \{ \phi(s_B), \min\{ \phi(t) - \phi(s) | (s,t) \in E(T) \ \mbox{and} \ s \in B \}\}.$$
\end{definition}

\begin{example}
In the tree in Example \ref{ex-tree}, the minimum width of the branches are $$\min\{d_1, d_{i+1}-d_i \ \mbox{for} \  1\leq i \leq 3\} \quad \mbox{and} \quad \min\{d_5, d_4-d_5\},$$ respectively. In Example \ref{ex-flag}, the minimum width of the branch is $$\min\{k_1, n-k_m, k_{i+1}-k_i \ \mbox{for} \ 1 \leq i \leq m-1\}.$$
\end{example}

\subsubsection{Results on finiteness of orbits} Our first theorem classifies the tree varieties that are homogeneous or have two orbits. 

\begin{theorem}\label{thm-homogeneous}
Let $F(T, \phi)$ be a tree variety.
\begin{enumerate}
\item The following are equivalent:
\begin{enumerate}
\item[(i)] The variety $F(T, \phi)$ is homogeneous.
\item[(ii)] The tree $T$ is a chain.
\item[(iii)] The variety $F(T, \phi)$ is a partial flag variety.
\end{enumerate}
\item The variety $F(T, \phi)$ has two $\PGL(n)$ orbits if and only if $T$ has exactly two branches each of length $1$  and one of the branches has minimum width equal to $1$.
\end{enumerate}
\end{theorem}
This theorem is closely related to the following result of Knop \cite{Knop} in Type A and generalizes Case (1). 

\begin{proposition}\cite{Knop}\label{cor-Freitag}
Let $X$ be a projective rational homogeneous variety. Then the complement of the diagonals in $X^m$ is homogeneous if and only if 
\begin{enumerate}
\item either $m=2$ and $X\cong\PP^n$ 
\item or $m=3$ and $X \cong \PP^1$.
\end{enumerate}
\end{proposition}

We next classify tree varieties that have finitely many orbits under the $\PGL(n)$ action, completely answering Question (1).

\begin{theorem}\label{thm-finiteorbits}
The tree variety $F(T, \phi)$ has finitely many $\PGL(n)$ orbits if and only if $(T, \phi)$ has at most 3 leaves and satisfies one of the following.
\begin{enumerate}
\item $T$ has at most 2 leaves.
\item $T$ has 3 leaves with the following possible branch lengths.
\begin{enumerate}
\item $(1,1,\ell)$ with $1 \leq \ell$,
\item $(1,2,\ell)$ with $2 \leq \ell \leq 4$,
\item $(1,2, \ell)$ with $5 \leq \ell$ provided that the minimum width of the branch of length 1 is 2 or the minimum width of the branch of length 2 is 1.
\item $(1, \ell_1, \ell_2)$ with $1 \leq \ell_1 \leq \ell_2$ provided that the minimum width of the branch of length 1 is 1.
\end{enumerate}
\end{enumerate}
\end{theorem}

As a special case, this theorem contains Magyar, Weyman and Zelevinsky's theorem classifying products of flag varieties with finitely many $\PGL(n)$ orbits \cite[Theorem 2.2]{MWZ}. 
One can also enumerate all the orbits in the cases described in Theorem \ref{thm-finiteorbits} using \cite[Theorem 2.9]{MWZ}.

\subsubsection{Results on density of orbits} The original motivation for introducing tree varieties was to classify products of partial flag varieties that have a dense $\PGL(n)$ orbit. 

\begin{example}
Any tree variety with finitely many $\PGL(n)$ orbits has a dense orbit. Hence, Theorem \ref{thm-finiteorbits} provides many examples of dense tree varieties. However, a tree variety may have infinitely many orbits, but still have a dense orbit. For instance, consider the tree consisting of $k \leq n+1$ vertices labeled $1$ each connecting to the root labeled $n$. 
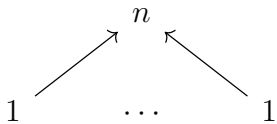
\[
\begin{tikzcd}
& n & \\
1\arrow[ru] & \cdots & 1\arrow[lu]  
\end{tikzcd}
\]
The corresponding variety is $k$ ordered points in $\PP^{n-1}$. This tree variety has a dense orbit since any linearly general $k$ points are equivalent under the $\PGL(n)$ action \cite[Exercise 1.6]{harris}, but has finitely many orbits only if $k \leq 3$.
\end{example}

If $\PGL(n)$ acts with dense orbit on a variety $X$, then we must have that $$\dim(\PGL (n)) = n^2 -1 \geq \dim(X).$$ This imposes strong dimension restrictions on tree varieties that can have a dense orbit. More generally, given any vertex $v \in T$, let $T^v$ be the subtree of $T$ consisting of the vertices that have a directed path terminating at $v$. Letting $\phi^v$ be the restriction of $\phi$ to $T^v$, we obtain a new labeled tree $(T^v, \phi^v)$ with root $v$.

\begin{lemma}\label{lem-stilltrivial}
If  $F(T, \phi)$ is dense, then for any vertex $v \in T$
$$ \sum_{(s,t) \in E(T^v)} \phi(s)(\phi(t) - \phi(s)) \leq \phi(v)^2-1.$$
\end{lemma}  

\begin{proof}
By comparing dimensions of stabilizers of general points and using Lemma \ref{lem-basiclem}, it follows that if $F(T, \phi)$ is dense, then $F(T^v, \phi^v)$ is dense. Hence, the dimension of the tree variety $F(T^v, \phi^v)$ has to be less than or equal to the dimension of $\PGL(\phi(v))$. By Theorem \ref{thm-irredanddim}, the dimension of the tree variety is given in the left-hand side of the inequality and the dimension of $\PGL(\phi(v))$ is given in the right-hand
side of the inequality. This proves the lemma.
\end{proof}
This motivates the following definition.

\begin{definition}
We call a tree variety $F(T, \phi)$ {\em dense} if $F(T, \phi)$ has a dense $\PGL(n)$ orbit. Otherwise, we say $F(T, \phi)$ is {\em sparse}. The tree variety $F(T, \phi)$ is {\em trivially sparse} if any vertex $v \in T$ violates the inequality in Lemma \ref{lem-stilltrivial}.
\end{definition}
Trivially sparse tree varieties are sparse for an easy to check reason. A tree variety may be sparse without being trivially sparse. The following is a generalization of \cite[Example 1.2]{CHZ}.
\begin{example}
Consider the labeled tree
\[
\begin{tikzcd}
1 \arrow[rd] \\
1\arrow[r] & m \arrow[r] & n \\
m-1 \arrow[ru] & m-1 \arrow[u]
\end{tikzcd}
\]
The tree variety associated to this tree is not trivially sparse when $m>2$, but it is sparse. Let $W$ be the span of the two $1$-dimensional subspaces. The vector space $W$ generically intersects each of  the $(m-1)$-dimensional subspaces in a $1$-dimensional subspace. The cross-ratio of the four $1$-dimensional subspaces in $W$ is an invariant of the orbits. This tree variety has dimension $4m-4+ m(n-m)$, which can be arbitrarily smaller than $n^2 -1$ as $n$ tends to infinity.
\end{example}

This example raises the following problem.

\begin{problem}
Classify dense tree varieties. 
\end{problem}

Already the following  special case seems to be  challenging. 
\begin{problem}
Classify dense tree varieties with three leaves.
\end{problem}

In Proposition \ref{prop-denseredprod}, we will show that classifying dense tree varieties with at most three leaves reduces to classifying dense products of three partial flag varieties. 
Popov \cite{Popov, Popov2} classified dense $(G/P)^n$ when $P$ is a maximal parabolic subgroup. Devyatov \cite{Devyatov}  has extended the  classification to non-maximal parabolic subgroups, except in type $A$. When $G/P$ is a type A partial flag variety, even the classification of dense $(G/P)^3$ is unknown. We will give several partial results towards this classification. Some of our results can be summarized in the following theorem.

\begin{theorem}\label{thm-introdense}
\begin{enumerate}
\item If there exists two indices $i \not= j$ such that $k_i + k_j =n$, then $F(k_1, \dots, k_r; n)^3$ is sparse (Corollary \ref{cor-addton}).
\item If $3k_r \leq n$, then $F(k_1, \dots, k_r; n)^3$ is dense (Lemma \ref{lem-easydense}).
\item If $2k_r \leq n$ and $2k_i \leq k_{i+1}$ for $2 \leq i \leq r-1$, then $F(k_1, \dots, k_r; n)^3$ is dense (Proposition \ref{prop-double}).
\end{enumerate}
\end{theorem}

Finally, we will classify the density of the triple self-product of two-step flag varieties.

\begin{theorem}[Proposition \ref{exprop-2step} and Theorem \ref{thm-2stepmain}]\label{thm-intro2step}
The product $F(k_1, k_2; n)^3$ is sparse if and only if $k_1 + k_2 = n$. The product $F(k_1, k_2; n)^3$ is trivially sparse if and only if $n$ is divisible by $3$, $k_1 = \frac{n}{3}$ and $k_2= \frac{2n}{3}$.
\end{theorem}

The density of the $\PGL(n)$ action on a product of flag varieties  has many applications. For the applications, we specialize our base field to $\CC$. Let $\lambda_1, \dots, \lambda_d$  be nonzero dominant characters of the maximal torus $T$ in the semi-simple group $G$. Then $(\lambda_1, \dots, \lambda_d)$  is called {\em primitive} if for every non-negative $d$-tuple of integers $(n_1, \dots, n_d)$, the Littlewood-Richardson coefficient $c^0_{n_1\lambda_1,\dots ,n_d \lambda_d} \leq 1$. Popov in \cite[Theorem 1]{Popov2} proves that if $G$ has an open orbit on $G/P_{\lambda_1} \times \cdots \times G/P_{\lambda_d}$, then $(\lambda_1, \dots, \lambda_d)$ is primitive. Hence, for vectors for which  $\PGL(n)$ acts with dense orbit, we get a strong bound on the Littlewood-Richardson coefficients. Similarly, the density has  geometric applications to enumerative geometry and genus zero Gromov-Witten invariants. We refer the reader to \cite{CHZ} for more details. 

Our approach to Theorems \ref{thm-finiteorbits} and \ref{thm-intro2step} is elementary. In order to show that a tree variety does not have a dense orbit, we explicitly construct a cross-ratio which has to be preserved by the $\PGL(n)$ action. For applications in other contexts, knowing the explicit cross-ratio which obstructs density is often useful.

\subsection*{Organization of the paper} In \S \ref{sec-prelim}, we recall the necessary background. In \S \ref{sec-feworbits}, we prove Theorem \ref{thm-homogeneous}. In \S \ref{sec-finite}, we classify tree varieties with finitely many orbits and prove Theorem \ref{thm-finiteorbits}. In \S \ref{sec-dense}, we study tree varieties with dense orbit.

\subsection*{Acknowledgements} We would like to thank Dave Anderson, James Freitag, Majid Hadian, J\'{a}nos Koll\'{a}r, Howard Nuer, Sybille Rosset, Geoffrey Smith and Dmitry Zakharov for helpful discussions regarding actions of $\PGL(n)$ on products of varieties. 

\section{Preliminaries}\label{sec-prelim}

In this section, we collect basic facts concerning tree-varieties and group actions.

\subsection{Tree-varieties} Let $(T, \phi)$ be a labeled tree with the root $r$. Recall that $V(T)$ and $E(T)$ denote the vertices and edges of $T$, respectively. Let $n = \phi(r)$ be the ambient dimension of the tree. The tree variety $F(T, \phi)$ parameterizes subspaces $(U_v)_{v \in V(T)}$ such that $\dim(U_v) = \phi(v)$ and $U_s \subset U_t$ whenever $(s,t) \in E(T)$. The tree variety $F(T, \phi))$ is naturally a closed algebraic subset of $\prod_{v \in V(T) \backslash \{r\}} G(\phi(v), n)$ given by imposing the incidence relations $U_s \subset U_t$ for every edge $(s,t) \in E(T)$.

\begin{theorem}\label{thm-irredanddim}
The tree variety $F(T, \phi)$ is a smooth, projective, irreducible variety of dimension $$\sum_{(s,t) \in E(T)} \phi(s) (\phi(t) - \phi(s)).$$
\end{theorem}
\begin{proof}
We prove the theorem by induction on the number of vertices in $T$. If $T$ has only two vertices and one edge $(s,t)$, then $F(T, \phi)$ parameterizes $\phi(s)$-dimensional subspaces of the $n$-dimensional vector space $W$. In this case, $F(T, \phi)$ is the Grassmannian $G(\phi(s), n)$ which is a smooth, irreducible, projective variety of dimension $\phi(s) (n - \phi(s))$. Since $n = \phi(t)$, the theorem is true in this case. 

By induction, assume that the theorem holds for trees with $m$ or fewer vertices. Let $T$ be a tree with $m+1$ vertices. Then $T$ has at least one leaf. Let $s_0$ be a leaf. Removing the leaf $s_0$ and the edge $(s_0,t_0)$  with source $s_0$, we obtain a tree $T'$ with $m$ vertices. The restriction of $\phi$ to $T'$ defines a function $\phi'$. Since $T'$ has $m$ vertices, by induction, the tree variety $F(T', \phi')$ is a smooth, irreducible, projective variety of dimension $\sum_{(s,t) \in E(T')} \phi(s) (\phi(t)-\phi(s))$. The variety $F(T, \phi)$ is obtained from $F(T', \phi')$ by choosing a $\phi(s_0)$-dimensional linear space in $U_{t_0}$. Hence, $F(T, \phi)$ is naturally a $G(\phi(s_0), \phi(t_0))$-bundle over $F(T', \phi')$. 
Consequently, $F(T, \phi)$ is a smooth, irreducible, projective variety with dimension $\phi(s_0) (\phi(t_0) - \phi(s_0)) + \dim (F(T', \phi'))$. The latter expression is precisely $\sum_{(s,t) \in E(T)} \phi(s) (\phi(t) - \phi(s)).$ The theorem follows by induction. 
\end{proof}

\subsubsection{Forgetful morphisms between tree varieties} Let $v$ be a vertex of $T$ different from the root. The vertex $v$ may be the target of more than one edge of $T$; however, $v$ is the source of a single edge $(v, t)$. Let $T_v$ be the tree obtained from $T$ by deleting $v$ and replacing every edge $(s, v)$ whose target is $v$ by $(s, t)$. Given a labeled tree $(T, \phi)$, we obtain a new labeled tree $(T_v, \phi_v)$, where $\phi_v$ is the restriction of $\phi$ to $V(T) \setminus v$. Then there is a natural forgetful morphism $$\pi_v: F(T, \phi) \to F(T_v, \phi_v)$$ that forgets the linear space $U_v$. This map is induced by the natural projection $$\prod_{w \in V(T) \backslash \{r\}} G(\phi(w), n) \to  \prod_{w \in V(T) \backslash \{v, r\}} G(\phi(w), n).$$

Given any vertex $v$ in $T$, there is a unique chain connecting $v$ to the root. Given a set of vertices $v_1, \dots, v_{\ell}$,  let $T_{v_1, \dots, v_{\ell}}$ denote the tree obtained from $T$ by deleting the vertices $v_1, \dots, v_{\ell}$ and replacing any edge $(s, v_i)$ with $(s, t_i)$, where $t_i$ is the first vertex in the chain connecting $s_i$ to the root which is not among $v_1, \dots, v_{\ell}$. Let $\phi_{v_1, \dots, v_{\ell}}$ be the restriction of $\phi$ to $V(T) \setminus \{v_1, \dots, v_{\ell}\}$.  Then there is a natural forgetful morphism $$\pi_{v_1, \dots, v_{\ell}}: F(T, \phi) \to F(T_{v_1, \dots, v_{\ell}}, \phi_{v_1, \dots, v_{\ell}})$$ that forgets the linear spaces $U_{v_1}, \dots, U_{v_{\ell}}$. This morphism is also induced by the corresponding natural projection 
$$\prod_{w \in V(T) \backslash \{r\}} G(\phi(w), n) \to  \prod_{w \in V(T) \backslash \{v_1, \dots, v_{\ell}, r\}} G(\phi(w), n).$$

\begin{proposition}\label{prop-forgetful}
Let $v$ be a vertex of $T$ different from the root. Let $s_1, \dots, s_j$ be the vertices of $T$ such that $(s_i, v)$ are edges in $T$. Then the forgetful morphism 
$$\pi_v : F(T, \phi) \to F(T_v, \phi_v)$$ is surjective if and only if $$\sum_{i=1}^j \phi(s_i) \leq \phi(v).$$ In particular,  if $v$ is the target of a unique edge, then $\pi_v$ is surjective and the fibers of $\pi_v$ are isomorphic to Grassmannians $G(\phi(v)-\phi(s_1), \phi(t) - \phi(s_1))$. 
\end{proposition}

\begin{proof}
A point $\{U_w\}_{w \in V(T_v)}$ of $F(T_v, \phi_v)$ is in the image of $\pi_v$ if and only if there is a linear space $U_v$ of dimension $\phi(v)$ contained in $U_t$ and containing $U_{s_i}$ for $1 \leq i \leq j$. If $\sum_{i=1}^j \phi(s_i) \leq \phi(v)$, one can always choose such a linear space $U_v$. Conversely, if $\sum_{i=1}^j \phi(s_i) > \phi(v)$, then $F(T_v, \phi_v)$ will contain points where the linear spaces $U_{s_i}$ span a vector space of dimension greater than $\phi(v)$. Hence, such a point cannot be in the image of $\pi_v$.
\end{proof}

\subsubsection{Constructing tree varieties inductively}
Given a labeled tree $(T, \phi)$ and a vertex $s \in V(T)$ different from the root $r$, there is a unique chain connecting $s$ to $r$. Define the distance function $d: V(T) \to \mathbb{N}$ by setting $d(s)$ to be the length of this chain for $s \not= r$ and set $d(r)=0$. 

Given a positive integer $m$, we can define {\em the truncation $(T_{\leq m}, \phi_{\leq m})$  of $(T, \phi)$ at distance $m$} as follows.  Let $T_{\leq m}$ be the tree obtained by deleting all the vertices $s$ of $T$ with $d(s) > m$ and deleting the edges that have these vertices as sources. Define  $\phi_{\leq m}$ by restricting $\phi$ to vertices $v$ with $d(v) \leq m$. If we delete the vertices of $T$ with  $d(v) < m$, we obtain a set of labeled trees $(T^1, \phi^1), \dots, (T^j, \phi^j)$ one for each vertex $v_i$ with $d(v_i) = m$. The vertex $v_i$ forms the root of the tree $T^i$ and $\phi^i$ is the restriction of $\phi$ to vertices at a distance at least $m$ that have $v_i$ in the chain connecting them to the root of $T$.

\begin{proposition}\label{prop-buildinductive}
For every positive integer $m$, there is a surjective forgetful morphism $$\pi_{\leq m}: F(T, \phi) \to F(T_{\leq m} , \phi_{\leq m})$$ and the fibers are isomorphic to $$F(T^1, \phi^1)\times \cdots \times  F(T^j, \phi^j).$$ In particular, a tree variety can be constructed inductively according to the distance function. 
\end{proposition}

\begin{proof}
The forgetful morphism $\pi_{\leq m}$ forgets all the vector spaces associated to vertices $v$ with $d(v) > m$. Let $v_1, \dots, v_j$ be the vertices with $d(v_i)=m$. Given a point in $F(T_{\leq m} , \phi_{\leq m})$, the fiber of $\pi_{\leq m}$ corresponds to choosing linear subspaces in $U_{v_i}$ according to the labeled tree $(T^i, \phi^i)$. The proposition follows.
\end{proof}

\subsection{Group actions}
\begin{lemma}\label{lem-basiclem}
Let $X$ be an irreducible projective variety with a $\PGL(n)$ action. Let  $x \in X$ be a closed point and let $\Stab(x)$ denote the stabilizer of $x$. Then the orbit of $x$ is dense in $X$ if and only if $$\dim(\Stab(x)) = n^2-1 - \dim(X).$$
\end{lemma}
\begin{proof}
Let $G$ be an algebraic group acting on an irreducible projective variety $X$. Then the orbit $Gx$ of $x$ under $G$ is open in its Zariski closure $\overline Gx$ by \cite[I.1.8]{Borel}. On the other hand, $Gx$ is isomorphic to $G/\Stab(x)$. Hence
$$\dim(\overline{Gx})= \dim(Gx) = \dim(G) - \dim(\Stab(x)).$$ Since $X$ is irreducible, the orbit $Gx$ is dense in $X$ if and only if $\dim(Gx) = \dim(X)$. Hence, the orbit of $x$ is dense if and only if $$\dim(\Stab(x))= \dim(G) - \dim(X).$$ The lemma follows by letting $G= \PGL(n)$ and noting that $\dim(\PGL(n))= n^2 -1$.  
\end{proof}

\begin{proposition}\label{prop-compare}
Let $\pi_v: F(T, \phi) \to F(T_v, \phi_v)$ be a surjective forgetful morphism. 
\begin{enumerate}
\item If $F(T, \phi)$ has a dense $\PGL(n)$ orbit, then $F(T_v, \phi_v)$ has a dense $\PGL(n)$ orbit. 
\item If $F(T, \phi)$ has finitely many $\PGL(n)$ orbits, then $F(T_v, \phi_v)$ has finitely many $\PGL(n)$ orbits. 
\end{enumerate}
\end{proposition}

\begin{proof}
The forgetful morphism $\pi_v$ is equivariant for the $\PGL(n)$ action. Hence, the image of an orbit is contained in an orbit. Let $O_T \subset F(T, \phi)$ be the dense orbit. Let $O \subset F(T_v, \phi_v)$ be the orbit containing $\pi_v(O_T)$. We then have $$\overline{O} \supset \pi_v(\overline{O_T}) = \pi_v(F(T, \phi))= F(T_v, \phi_v).$$ This proves part (1).

Suppose $F(T, \phi) = \sqcup_{i=1}^j O_i$ is a union of finitely many $\PGL(n)$ orbits. The image $\pi_v(O_i)$ is contained in an orbit $O_i'$. Then $$F(T_v, \phi_v) = \pi_v(F(T, \phi)) =  \sqcup_{i=1}^j \pi_v(O_i)  \subset \cup_{i=1}^j O_i'.$$ Hence, $F(T_v, \phi_v)$ has finitely many orbits. Of course, some of the orbits $O_i'$ may coincide. This concludes the proof of the proposition.
\end{proof}

\begin{proposition}\label{prop-complement}
The action of $\PGL(n)$  on $\prod_{i=1}^\ell F(k_{i, 1}, \dots, k_{i, j_i}, n)$ has finitely many orbits (respectively, a dense orbit) if and only if the action of $\PGL(n)$ on  $\prod_{i=1}^\ell F(n-k_{i, j_1}, \dots, n-k_{i, 1}, n)$ has finitely many orbits (respectively, a dense orbit).
\end{proposition}

\begin{proof}
Let $W^*$ be the dual of the ambient vector space $W$ with the dual $\PGL(n)$ action. Taking quotient spaces and passing to the dual defines an isomorphism between $\prod_{i=1}^\ell F(k_{i, 1}, \dots, k_{i, j_i}, n)$ and $\prod_{i=1}^\ell F(n-k_{i, j_1}, \dots, n-k_{i, 1}, n)$ which respects the $\PGL(n)$ action. The proposition follows. 
\end{proof}

The action of $\PP GL(n)$ on products of Grassmannians has been studied in detail in \cite{CHZ}. We recall the following theorem for the reader's convenience (see also \cite{SW}).

\begin{theorem}\cite[Theorem 5.1]{CHZ}
Let $m \leq 4$ and let  $X=\prod_{i=1}^m G(k_i, n)$ be a product of $m$-Grassmannians. Then the $\PP GL (n)$  action on $X$ is sparse if and only if $m=4$ and $\sum_{i=1}^4 k_i = 2n$.
\end{theorem} 

The first case of the theorem is the action of $\PP GL(2)$ on an ordered set of $4$ distinct points $(z_1, z_2, z_3, z_4)$. In this case, the action is trivially sparse and the invariants are fully understood. There is a unique element $g$ of $\PP GL(2)$ taking the first three to $0, \infty$ and  $1$, respectively. The image of the fourth point $g(z_4)$  is called the cross-ratio of the four points. Four distinct ordered points are projectively equivalent if and only if their cross-ratio is the same.

\section{Tree varieties with few $\PGL(n)$ orbits}\label{sec-feworbits}

In this section, prove Theorem \ref{thm-homogeneous} and Proposition \ref{cor-Freitag}. 

\begin{proof}[Proof of Theorem \ref{thm-homogeneous}] 
We first prove (1). If the tree $T$ is a chain, then, by Example \ref{ex-flag}, $F(T, \phi)$ is a partial flag variety.  Hence, $F(T, \phi)$ is homogeneous under the $\PGL(n)$ action. This shows (ii) implies (iii) implies (i). To conclude the proof of (1), we need to show that (i) implies (ii).

Suppose that the tree $T$ is not a chain. Then there exists a vertex $v$ such that there are at least two vertices $v_1$ and $v_2$ so that $(v_1,v)$ and $(v_2,v)$ are edges in $T$. Let $t$ be such a vertex with the smallest distance from the root.   Let 
 $$\phi(v_1) = d_1, \quad \phi(v_2)=d_2 \quad \mbox{and}  \quad \phi(t)=d_t.$$  Without loss of generality, we may assume that $d_1 \leq d_2 < d_t$. Among the linear spaces parameterized by the  tree variety $F(T, \phi)$, there are  three linear spaces $U_1$, $U_2$, $U_t$ of dimensions $d_1$, $d_2$ and $d_t$, respectively, corresponding to these three vertices.  By Proposition \ref{prop-buildinductive}, we may construct $F(T, \phi)$ inductively starting at the root. Once we have chosen $U_t$, $U_1$ and $U_2$ are arbitrary linear subspaces of $U_t$ of dimensions $d_1$ and $d_2$, respectively. Hence, $$\max(0, d_1 + d_2 - d_t) \leq \dim(U_1 \cap U_2)  \leq d_1$$ and every possible value in this range can occur. Since $d_1 > 0$ and $d_2<d_t$, there are at least two possible values. Since the dimension $\dim(U_1 \cap U_2)$ is an invariant of the $\PGL(n)$ action, $F(T, \phi)$ cannot be homogeneous. We conclude that (i) implies (ii).  
 
Next suppose $F(T, \phi)$ has two orbits under the $\PGL(n)$ action. Then $T$ is not a chain and we do have vertices $v_1, v_2$ and $t$ as in the previous paragraph. In this case, the range  $$\max(0, d_1 + d_2 - d_t) \leq \dim(U_1 \cap U_2)  \leq d_1$$ must have 2 possible values. This can only happen if $d_1=1$ or $d_2 = d_t -1$. To conclude the proof of (2), we need to show that $v_1$ and $v_2$ are leaves of $T$ and there are no other edges with $t$ as the target. 

Suppose there is a third vertex $v_3$ such that $(v_3, t) \in E(T)$. Without loss of generality, we may assume that $\phi(v_3) = d_3 \geq d_2$.  The corresponding vector space $U_3 \subset U_t$ may be chosen freely in $U_t$. If $U_1 \subset U_2$, then $U_3$ may or may not contain $U_2$.  Further, if $U_3$ does not contain $U_2$, it may or may not contain $U_1$. If $U_1 \not\subset U_2$, then $U_3$ may or may not contain either $U_1$ or $U_2$. We conclude that there are at least 6 orbits. Since $F(T, \phi)$ has only two $\PGL(n)$ orbits, there cannot be a third vertex $v_3$ such that $(v_3, t) \in E(T)$. 

If there is an edge $(v_3, v_1) \in E(T)$, then $U_3 \subset U_1$ can be chosen freely. Then there are at least three possibilities:
\begin{enumerate}
 \item  $U_3 \subset U_1 \subset U_2$,  
 \item  $U_3 \subset U_1$ and neither are subsets of $U_2$, or 
 \item $U_3 \subset U_1 \cap U_2$, but $U_1$ is not a subset of $U_2$.
 \end{enumerate}
 Hence, there are at least 3 orbits. A similar argument applies if there is an edge $(v_3, v_2)$. We conclude that if $F(T, \phi)$ has two orbits under the $\PGL(n)$ action, then the tree looks like

\[
\begin{tikzcd}
d_1 \arrow[r] & d_3 \arrow[r] & d_4  \arrow[r] & \cdots  \arrow[r]& n \\
& d_2  \arrow[u]& & &   &
\end{tikzcd}
\]
and assuming that $d_1 \leq d_2 \leq d_3$, either $d_1=1$ or $d_2 = d_3 -1$. Hence, the tree has exactly two branches each of length one and one of the branches has minimum width equal to $1$.

Conversely, the $\PGL(n)$ action has two orbits on such a tree variety. The orbits are determined by whether $U_1 \subset U_2$ or $U_1 \not\subset U_2$. If $U_1 \subset U_2$, then 
the corresponding orbit is a partial flag variety, hence homogeneous. If $U_1 \not\subset U_2$, then $U_2 \subset U_3 \subset \cdots \subset W$ is a partial flag and we may choose a basis for $W$ so that $U_i$ is the span of $e_i$ for $1 \leq i \leq d_i$. If $\dim(U_1)=d_1=1$, then we may further require that $e_{d_3}$ is a basis for $U_1$. If $d_1 >1$ and $d_3-d_2=1$, we may require that  $e_i$ for $d_2-d_1+2\leq i \leq d_2$ is a basis for $U_1 \cap U_2$ and $U_1$ is spanned by $e_i$ for $d_2-d_1+2\leq i \leq d_3$. Hence, this locus also forms a single orbit. We conclude that these tree varieties have exactly two $\PGL(n)$ orbits. This proves part (2) of  Theorem \ref{thm-homogeneous}.
\end{proof}

For completeness, we sketch a simple proof of Knop's result \cite{Knop}. 

\begin{proof}[Proof of Proposition \ref{cor-Freitag}]
By Poincar\'e duality, a  one-dimensional Schubert variety in a rational homogeneous variety $X= G/P$ is a line in the minimal embedding of $X$. Suppose that $m>1$ and that the complement of the diagonals in $X^m$ is homogeneous. Observe that this implies that the complement of the diagonals in $X^{l}$ is homogeneous for all $l \leq m$. In particular, the complement of the diagonals in $X^2$ is homogeneous. Pick a line $L$ on $X$ and let  $p$ and $q$ be distinct points on $L$. Since $G$ is acting transitively on pairs of points, there must be a line between any two distinct points on $X$. We conclude that $X= \PP^n$ for some $n$. If $n>1$, we can take three distinct collinear points and three distinct non-collinear points on $\PP^n$ to see that the complement of the diagonals in  $(\PP^n)^3$ is not homogeneous. When $n=1$,  $\PGL(2)$ acts transitively  on triples of ordered, distinct points on $\PP^1$. Since $3 = \dim \PGL(2) < \dim((\PP^1)^4) = 4$, $\PGL(2)$ does not have a dense orbit on $(\PP^1)^4$. This concludes the proof of the proposition. 
\end{proof}

\section{Tree varieties with finitely many $\PGL(n)$ orbits}\label{sec-finite}
In this section, we classify tree varieties with finitely many $\PGL(n)$ orbits and prove Theorem \ref{thm-finiteorbits}. It is possible to reduce the proof of this theorem to the classification of flag varieties of finite type by  Magyar, Weyman and Zelevinsky \cite{MWZ} and obtain a relatively short proof. However, we prefer to give an elementary proof.

\begin{proof}[Proof of Theorem \ref{thm-finiteorbits}]
We will classify  $F(T, \phi)$ that have finitely many $\PGL(n)$ orbits.  The classification is somewhat involved, so we will break it into smaller steps. From now on suppose that  $F(T, \phi)$ has  finitely many $\PGL(n)$ orbits. We begin by showing that $T$ must have at most 3 leaves. We will then show that if $T$ has at most 2 leaves, then $F(T, \phi)$ has finitely many orbits. The hardest part is to classify trees with 3 leaves for which $F(T, \phi)$ has finitely many orbits. To show that $F(T, \phi)$ does not have finitely many $\PGL(n)$ orbits, we construct a cross-ratio that needs to be preserved but can take arbitrary values.

\subsection*{Step 1: $T$ has at most 3 leaves}  We first show that if $F(T, \phi)$ has finitely many $\PGL(n)$ orbits, then $T$ has at most 3 leaves. Suppose that $T$ has  4 leaves $s_1, \dots, s_4$. For each leaf $s_i$, let $(s_i, t_i)$ be the edge with source $s_i$. Some of the vertices $t_i$ may coincide. Fix a full flag $$F_1 \subset \cdots \subset F_n=W,$$ where $F_k$ has dimension $k$. For each vertex $v \in V(T)$ different from the leaves $s_1, \dots, s_4$, let $U_v = F_{\phi(v)}$. Then the 4 linear spaces $U_i$ corresponding to the leaves $s_i$ can be chosen freely subject to the condition that $U_i \subset F_{\phi(t_{i})}$. 

Let $d_i$ denote the dimension of $U_i$ and without loss of generality assume that $d_1 \leq d_2 \leq d_3 \leq d_4$. Choose $U_2 \subset \bigcap_{i=2}^4 F_{\phi(t_i)}$ so that $\dim(U_2 \cap F_{\phi(t_1)}) \geq d_1 -1$. Since $\min(\phi(t_3), \phi(t_4)) > d_3$ and $\min(\phi(t_1), \phi(t_2)) > d_1$, this is possible. Choose $U_1 \subset \bigcap_{i=1}^4 F_{\phi(t_i)}$ such that $U_1 \cap U_2 = \Lambda_1$ with $\dim(\Lambda_1) = d_1 -1$ and the span of $U_1$ and $U_2$ is $\Lambda_2$ with $\dim(\Lambda_2) = d_2 + 1$. Let $U_3 \subset F_{\phi(t_{3})}$ be a linear space of dimension $d_3$ containing $\Lambda_1$ and intersecting $\Lambda_2$ in a linear space of dimension $d_2$ not containing $U_1$ or $U_2$. Set $\Lambda_3 = U_2 \cap U_3$. Then $\dim(\Lambda_3) = d_2 -1$ and $\Lambda_1 \subset \Lambda_3$ by construction. Finally, choose  $U_4 \subset F_{\phi(t_{4})}$ of dimension $d_4$ containing $\Lambda_3$, intersecting $\Lambda_2$ in a linear space  of dimension $d_2$ and not containing $U_1, U_2, U_3$. For each $i$, let $Z_i$ denote the span of $\Lambda_3$ and $U_i \cap \Lambda_2$. Observe that $\dim(Z_i) = d_2$ and $\Lambda_3 \subset Z_i$.
 
The action of $\PGL(n)$ respects spans and intersections, consequently $\PGL(n)$ acts on the vector spaces $Z_i$. The linear spaces $Z_i$ of dimension $d_2$ are all contained in $\Lambda_2$ of dimension $d_2 +1$ and contain $\Lambda_3$ of dimension $d_2-1$. Hence, they determine 4 points in a pencil of linear spaces. (Equivalently, $Z_i/\Lambda_3$ are 4  two-dimensional linear subspaces of  $\Lambda_2/\Lambda_3$.) The cross-ratio of these 4 points on $\PP^1$ is an invariant of the $\PGL(n)$ action. Since we can choose $Z_2$ arbitrarily subject to the condition that it contains $\Lambda_3$, is contained in $\Lambda_2$ and does not contain $U_1, Z_3, Z_4$, any cross-ratio is possible. Hence, $\PGL(n)$ cannot have finitely many orbits as soon as $T$ contains at least $4$ leaves and the base field is infinite.

\subsection*{Step 2: If $T$ has at most 2 leaves, then $F(T, \phi)$ has finitely many orbits.}  If $T$ has only one leaf, then $T$ is a chain. In this case $F(T, \phi)$ is a partial flag variety and is homogeneous under the $\PGL(n)$ action (see Example \ref{ex-flag} and Theorem \ref{thm-homogeneous} (1)). 

Next, suppose $T$ has two leaves, then $T$ has to be of the form
\[
\begin{tikzcd}
d_1 \arrow[r] & d_2 \arrow[r] & \cdots  \arrow[r] & d_j \arrow [r] & \cdots \arrow[r] & n  \\ 
& e_1 \arrow[r] & \cdots \arrow[r] & e_l  \arrow[u] & & &   
\end{tikzcd}
\]
Let $(T', \phi')$ be the following labeled tree
\[
\begin{tikzcd}
1 \arrow[r] & 2 \arrow[r] & \cdots  \arrow[r] & d_j-1 \arrow[r] & d_j \arrow [r]  & d_{j}+1  \arrow[r] & \cdots  \arrow[r]& n-1 \arrow[r] & n  \\ 
& & e_1 \arrow[r] & \cdots \arrow[r] & e_l  \arrow[u] & & &   
\end{tikzcd}
\]
By repeated applications of Proposition \ref{prop-forgetful}, there is a surjective  morphism $\pi: F(T', \phi') \rightarrow F(T, \phi)$. By Proposition \ref{prop-compare}, it suffices to show that $F(T', \phi')$ has finitely many $\PGL(n)$ orbits. Fix a basis $e_1, \dots, e_n$ for $W$. Since $\PGL(n)$ acts transitively on full flags, after applying an element of $\PGL(n)$ we may assume that the vector spaces $U_i$ corresponding to the top chain in $T'$ are given by the span of $e_j$ for $1 \leq j \leq i$.  The stabilizer of the full flag is the Borel subgroup of upper triangular matrices. Next, we need to choose a partial flag $W_{e_1} \subset \cdots \subset  W_{e_l} \subset U_{d_j}$. The orbits of the Borel group are precisely the Schubert cells of the partial flag variety $F(e_1, \dots, e_l, d_j)$. A Schubert cell is determined by specifying an element of the symmetric group $\mathfrak{S}_{d_j}$ with at most $l$-descents at $e_1, \dots, e_l$. There are finitely many Schubert cells, hence $\PGL(n)$ acts with finitely many orbits on $F(T', \phi')$. More generally, the orbits in the original tree variety $F(T, \phi)$ are determined by  the dimensions of intersections of $W_{e_\alpha}$ and $U_{d_\beta}$ for every $1 \leq \alpha \leq l$ and $1 \leq \beta \leq j-1$.
\smallskip

This concludes the discussion of trees with greater than or equal to 4 or less than or equal to 2 leaves. From now on, we assume that $T$ has three leaves.

\subsection*{Step 3: One of the branches has length at most 1} We next show that if $F(T, \phi)$ has finitely many $\PGL(n)$ orbits, then $T$ cannot have three branches where each branch has  at least 2 vertices. By repeatedly applying Proposition \ref{prop-forgetful}, we may forget all but two of the vector spaces comprising any branch of $T$ to obtain a tree $T'$ and a surjective morphism $\pi: F(T, \phi) \rightarrow F(T', \phi')$. By Proposition \ref{prop-compare}, it suffices to show that $F(T', \phi')$ does not have finitely many $\PGL(n)$ orbits. Hence we may assume that $T$ looks like 
\[
\begin{tikzcd}
s_1 \arrow[r] &t_1  \arrow[r] & v  \arrow[r] & \cdots \arrow[r] & v' \arrow[r] & \cdots \arrow[r] & n, \\
s_2 \arrow[r] & t_2   \arrow[ru]&  s_3 \arrow[r]& t_3 \arrow[ru]   & & 
\end{tikzcd}
\]
where the vertices labeled $v$ and $v'$ may coincide. 

Fix a flag $F_{\bullet}$ and choose all the vector spaces except for the 6 labeled $s_i, t_i$ for $1 \leq i \leq 3$ from this fixed flag. Now we can choose the 6 vector spaces  $U_{s_i} \subset U_{t_i}$ for $1 \leq i \leq 3$ only subject to the condition that $U_{t_i}$ is contained in the flag element $F_{\phi(v)}$ for $i=1,2$ and $U_{t_3}$ is contained in  the flag element $F_{\phi(v')}$.  In particular, we may assume that $U_{s_3}, U_{t_3}$ intersect $F_{\phi(v)}$ in distinct and proper subspaces. The proof of this step follows from the following claim. 

\begin{claim}\label{claim-222}
Let  $U_i \subset V_i$ for $1 \leq i \leq 3$, be $6$ proper subspaces of a vector space $W$. Then $\PGL(W)$ does not have finitely many orbits on such 6-tuples of vector spaces.
\end{claim}

\begin{proof}[Proof of Claim \ref{claim-222}]
Let the small letter denote the dimension of the vector space denoted by the capital letter.  The prototypical case is when $u_i=1$ and $v_i=2$ for all $i$. In that case, take $U_i \subset V_i$ to be general flags contained in a three dimensional subspace $Y$. Then $V_2$ and $V_3$ intersect $V_1$ in one-dimensional subspaces $Q_2$ and $Q_3$. Similarly, the span of $U_2$ and $U_3$ intersects $V_1$ in a one-dimensional subspace $Q_4$, distinct from $Q_2$ and $Q_3$. Then the cross-ratio of $U_1, Q_2, Q_3, Q_4$ is an invariant of the $\PGL(n)$ action. By varying $U_1$, we see that the cross-ratio takes arbitrary values. Hence, there are  infinitely many $\PGL(n)$ orbits if the base field is infinite.

The general case is similar.  Without loss of generality, assume that $v_1 \leq v_2 \leq v_3$. Fix a linear space $V_1$ of dimension $v_1$ and choose $V_2$ and $V_3$ to be general among linear spaces of dimensions $v_2$ and $v_3$, respectively, such that they intersect $V_1$ in subspaces $Q_2$ and $Q_3$ of dimension $v_1 -1$ and span a linear space $Y$ of dimension $v_3 + 1$. Observe that $Y$ contains $V_1$. Let $\Lambda = \cap_{i=1}^3 V_i$, which is a linear space of dimension $v_1 -2$. Let $U_1$ be a general linear subspace of $V_1$ which intersects $\Lambda$ in a subspace of dimension $\min (u_1 -1, v_1-2)$.  Let $Q_1 \subset V_1$ be the $(v_1-1)$-dimensional space spanned by $U_1$ and $\Lambda$.  Let $Q_4$ be a $(v_1-1)$-dimensional subspace of $V_1$ containing $\Lambda$ and distinct from $Q_1$, $Q_2$ and $Q_3$. Let $\gamma \in Q_4$ be a vector not contained in $\Lambda$. Let $\alpha_2$ be a vector in $V_2$ not contained in $V_1$ or $V_3$. Then there is a unique vector $\alpha_3 \in V_3$ such that $\gamma$ is a linear combination of $\alpha_2$ and $\alpha_3$. Let $Y'$ be a general codimension one linear subspace of $Y$ containing $Q_4$ and $\alpha_2$ (and hence also $\alpha_3$). 
For $i=2,3$, let $U_i$  be general linear subspaces of $Y' \cap V_i$ of dimension $u_i$ containing $\alpha_i$. Then the span of $U_2$, $U_3$ and $\Lambda$  intersects $V_1$ in $Q_4$. The linear spaces $Q_1, Q_2, Q_3, Q_4$ are $(v_1-1)$-dimensional subspaces of $V_1$ and contain $\Lambda$. Hence, they form a pencil. The cross-ratio of these four vector spaces in the pencil is an invariant of the $\PGL(n)$ action and can take arbitrary values. We conclude that $\PGL(n)$ cannot have finitely many orbits.
\end{proof}

We thus conclude that if $F(T, \phi)$ has finitely many $\PGL(n)$ orbits, then $T$ has at most three leaves and one of the branches has length at most one.

\subsection*{Step 4: Bounding the length of the other branches} Assume that $T$ has three leaves. We now show that  if the minimum width of the branch with one leaf is not one, then the length of another branch has to be at most 2. Suppose that $T$ has two branches $u_1 \to \cdots \to u_p$ and $v_1 \to \cdots \to v_q$ with $p, q \geq 3$ and another branch of length one with leaf $s_1$ and minimum width greater than 1. We would like to show that $F(T, \phi)$ does not have finitely many $\PGL (n)$ orbits.  By repeated applications of Propositions  \ref{prop-forgetful} and \ref{prop-compare}, it suffices to assume that $p=q=3$. Then the tree $T$ has one of the following forms: 
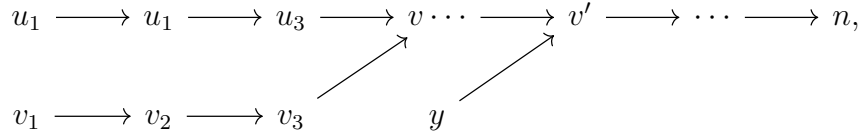
\[
\begin{tikzcd}
u_1 \arrow[r] &u_1  \arrow[r] & u_3  \arrow[r] & v \arrow[r] \cdots & v' \arrow[r] & \cdots \arrow[r] & n, \\
v_1 \arrow[r] & v_2   \arrow[r]& v_3 \arrow[ru]& y \arrow[ru]   & & 
\end{tikzcd}
\]
or

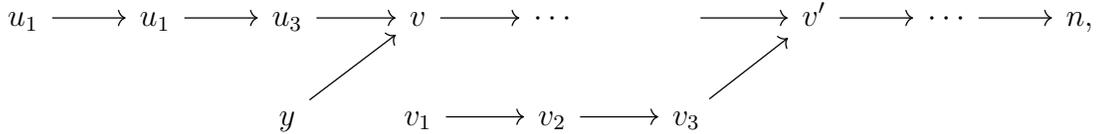
\[
\begin{tikzcd}
u_1 \arrow[r] &u_1  \arrow[r] & u_3  \arrow[r] & v \arrow[r] & \cdots & \arrow[r] & v' \arrow[r] & \cdots \arrow[r] & n, \\
& &  y \arrow[ru] &  v_1 \arrow[r] & v_2   \arrow[r]& v_3 \arrow[ru]&   & & 
\end{tikzcd}
\]

Let all the linear spaces corresponding to vertices other than those labeled by $u_i$, $v_i$ or $y$ be from a fixed flag. In the first case, let the vector space $Y$ of dimension $y$ intersect the vector space corresponding to the vertex $v$ in a subspace of dimension between $2$ and $v-2$. This is possible since the minimum width of that branch is greater than 1. In the second case, let the flag $V_1 \subset V_2 \subset V_3$ intersect the vector space corresponding to $v$ in a nontrivial three-step flag. Thus, the proof of this step reduces to the following claim.

\begin{claim}\label{claim-233}
Let $U_1 \subset U_2 \subset U_3$ and $V_1 \subset V_2 \subset V_3$ be proper linear spaces of a vector space $W$ of dimension $n$. Let $Y$ be a subspace of $W$ of dimension not equal to $1$ or $n-1$. Then $\PGL(W)$ cannot have finitely many orbits on the configurations of the $7$ vector spaces $U_i, V_i, Y$. 
\end{claim}
\begin{proof}[Proof of Claim \ref{claim-233}]
The prototypical case is when $u_1=v_1=1$, $u_2=v_2=2$ and $u_3=v_3=3$ and $y=2$. We may assume that these are general vector spaces in a 4-dimensional vector space. The span of $U_1$ and $V_2$ intersects $Y$ in a one-dimensional linear space $Q_1$. Similarly, the span of $U_2$ and $V_1$ intersects in a one-dimensional linear space $Q_2$. Finally, $U_3$ and $V_3$ intersect $Y$  in a one-dimensional linear spaces $Q_3$ and $Q_4$, respectively, The $\PGL(n)$ action needs to preserve the cross-ratio of the subspaces $Q_i$ in $Y$. Since the cross-ratio can take an arbitrary value, $\PGL(n)$ cannot act with finitely many orbits.

The general case is similar. First, we make two initial reductions to simplify the argument. Choose $V_2$ so that it intersects $U_2$ in a subspace of dimension $u_2 -2$. Let $\Omega$ be the span of $V_2$ and $U_2$, which has dimension $v_2 + 2$. We can then choose $U_3$, $V_3$ and $Y$ to intersect $\Omega$ in dimensions $\max(u_2 +1, u_3 + v_2 + 2-n)$, $v_2 + 1$ and $\max(2, y + v_2 + 2 -n)$, respectively. Since $\Omega$ is canonically determined as the span of $U_2$ and $V_2$, it suffices to show that $\PGL(n)$ has infinitely many orbits in the corresponding configuration in $\Omega$. We may thus assume that $v_2 = n-2$ and $v_3 = n-1$. 

Next, assume that $V_1$ and $U_3$ intersect in a linear space of dimension $\min(v_1-1, u_3-3)$. Let $\Sigma$ be the span of $V_1$ and $U_3$, which is a linear space of dimension $\sigma=\max(v_1+3, u_3+1)$. We can choose  $V_2, V_3$ and $Y$ to be linear spaces that intersect $\Sigma$ in dimension $\sigma-2$, $\sigma-1$ and $\max(2, y+\sigma-n)$, respectively. Since $\Sigma$ is canonically determined as the span of $U_3$ and $V_1$, it suffices to show that $\PGL(n)$ has infinitely many orbits in the corresponding configuration in $\Sigma$. We thus reduce to the case when either $(v_1, v_2, v_3) = (n-3, n-2, n-1)$ or $v_2=n-2, v_3=n-1$ and $u_3 = n-1$.

First, assume  $(v_1, v_2, v_3) = (n-3, n-2, n-1)$ and fix $V_1 \subset V_2$ of dimension $v_1$ and $v_2$. Let $U_1$ be a linear space of dimension $u_1$ general among linear spaces  which intersect $V_2$ in a subspace of dimension $u_1 -1$. Observe that the span of $U_1$ and $V_2$ is a linear space $Q_1'$ of dimension $n-1$. Let $U_2$ be a linear space of dimension $u_2$ general among linear spaces that contain $U_1$ and intersect $V_1$ in a subspace of dimension $u_2-2$. Observe that the span of $U_2$ and $V_1$ is a linear space $Q_2'$ of dimension $n-1$ distinct from $Q_1'$. Choose $V_3$ to be a general vector space of dimension $v_3$ containing $V_2$. Observe that it is distinct from $Q_1'$ and $Q_2'$. Let $\Lambda' = Q_1' \cap Q_2' \cap V_3$, which is a linear space of dimension $n-3$. Let $Y$ be a linear space of dimension $y$ general among those that intersect $\Lambda'$ in a subspace $\Lambda$ of dimension $y-2$. Set $Q_1 = Q_1' \cap Y$, $Q_2 = Q_2' \cap Y$ and $Q_3 = V_3 \cap Y$. Note that $Q_1, Q_2, Q_3$ are subspaces of $Y$ of dimension $y-1$ that contain $\Lambda$. Finally, pick a vector $v$  in $Y$ not contained in $Q_1, Q_2$ or  $Q_3$. Let $U_3$ be a linear space of dimension $u_3$ general among those that contain $v$, $U_2$ and intersect $Y$ along the span of $v$ and $\Lambda$. Let $Q_4$ be the span of $v$ and $\Lambda$, which is the span of $U_3 \cap Y$ with $\Lambda$. Then $Q_1, \dots, Q_4$ are four points of a pencil of hyperplanes in $Y$ and their cross-ratio is an invariant of the $\PGL(n)$ action. By varying $v$, we can get any cross-ratio. We conclude that $\PGL(n)$ cannot act with finitely many orbits.

Finally, assume that $v_2 = n-2, v_3 = n-1$ and $u_3 = n-1$. Fix $V_2 \subset V_3$ to be linear spaces of dimension $n-2$ and $n-1$, respectively. Let $U_1$ be a linear space of dimension $u_1$ general among those that intersect $V_2$ in a subspace of dimension $u_1 -1$. Then $U_1$ and $V_2$ span a linear space $Q_1'$ of dimension $n-1$, distinct from $V_3$. Let $U_3$ be a  linear space of dimension $n-1$ general among those that contain $U_1$. Let $\Lambda' = Q_1' \cap V_3 \cap U_3$, which is a linear space of $n-3$ contained in $V_2$. Let $Y$ be a linear space of dimension $y$ general among those which intersect $\Lambda'$ in a subspace $\Lambda$ of dimension $y-2$. Set $Q_1 = Q_1' \cap Y$, $Q_2 = V_3 \cap Y$ and $Q_3 = U_3 \cap Y$. Take a hyperplane general $H$ among those containing $U_1$ and $\Lambda$. Observe that this is possible since $Q_1'$ and $U_3$ both contain $U_1$. Let $V_1$ be $H \cap V_2$ and let $U_2$ be a general linear space containing $U_1$ and contained in $U_3 \cap H$. Then the span of $V_1$ and $U_2$ is $H$. Setting $Q_4 = H \cap Y$, the cross-ratio of $Q_1, \dots, Q_4$ in the pencil they span in $Y$ is an invariant of the $\PGL(n)$ action. We conclude that $\PGL(n)$ cannot act with finitely many orbits.
\end{proof}

\subsection*{Step 5: Bounding the length of the third branch}
Finally, we show that if $T$ is a tree with three leaves with branch lengths $1, 2$ and $\ell$ and the minimum width of the branch of length 1 is at least 3 or the minimum width of the branch of length 2 is at least 2, then $\ell \leq 4$. By repeated applications of Propositions \ref{prop-forgetful} and \ref{prop-compare}, it suffices to study the case when $\ell = 5$. By reductions similar to the previous cases, this step follows from the following claim. 

\begin{claim}\label{claim-12l}
Let $A$, $B_1 \subset B_2$ and $C_1 \subset \cdots \subset C_5$ be distinct, proper subspaces of a vector space $W$ of dimension $n$. Assume that $\dim(A) \not= 1, 2, n-2$ or $n-1$. Assume that $\dim(B_1) \not=1$, $\dim (B_2) \not= n-1$ and $\dim(B_2) \not= \dim(B_1) +1$. Then $\PGL(n)$ does not act on configurations of $8$ subspaces $A, B_i, C_j$ with finitely many orbits.
\end{claim}

\begin{proof}[Proof of Claim \ref{claim-12l}]
The prototypical example is when $n=6$, $\dim(A)=3$, $\dim (B_i) =2i$ and $\dim(C_i) = i$. Let all these linear spaces be general. For $i=1,2$, let $D_i$ be the span of $A$ and $C_i$. Let $U_1= A \cap B_2$, and for $i=2,3$, let $U_i = D_{i-1} \cap B_2$. We thus get three linear spaces $U_1\subset U_2 \subset U_3$, where $\dim(U_i) =i$ in the 4-dimensional linear space $B_2$.  For $i\geq 3$, let $V_{i-3} = B_2 \cap C_i$. We get another three linear spaces $V_1 \subset V_2 \subset V_3$ where $\dim(V_i) = i$ in $B_2$. Finally, let $Y = B_1$. It is easy to see that these are general linear spaces. We thus reduce to the configuration in Claim \ref{claim-233}, hence $\PGL(n)$ cannot act with finitely many orbits.

The general case is similar. Fix $B_1 \subset B_2$ to be general linear spaces of dimensions $b_1$ and $b_2$, respectively. Choose $C_3 \subset C_4 \subset C_5$ such that  $\dim (C_i \cap B_2) = \max(i-2, c_i + b_2 -n)$. Set $V_i = C_{i+2} \cap B_2$. Choose $C_1 \subset C_2$ subsets of $C_3$ and $A$ so that they satisfy the following conditions:
\begin{enumerate}
\item We have $\dim(A \cap B_2) = \max(1, a+b_2-n)$,
\item Let $D_i$ denote the span of $A$ and $C_i$. Then for $i=1,2$,  $\max(i+1, a+b_2 -n+i) \leq \dim(D_i \cap B_2)$
\end{enumerate}
Set $U_1 = A \cap B_2$ and for $i=2,3$ set $U_i = D_{i+1} \cap B_2$. Finally set $B_1 = Y$.  We obtain the configuration in Claim \ref{claim-233} and such configurations cannot have finitely many $\PGL(n)$ orbits.
\end{proof}

Now we are ready to complete the proof that if $F(T, \phi)$ has finitely many $\PGL(n)$ orbits, then $(T, \phi)$ must be one of the labeled trees listed in Theorem \ref{thm-finiteorbits}. If $T$ has at most two leaves, then, by Step 2, $F(T, \phi)$ has finitely many orbits. This corresponds to Case (1).  If $T$ has at least 4 leaves, then, by Step 1, $F(T, \phi)$ cannot have finitely many orbits. 

We now consider the case when $T$ has exactly three leaves.   By Step 3, one of the branches must have length 1. If the minimum width of the branch with length 1 is not 1, then Step 4 guarantees that another branch has length at most 2. Hence, if $F(T, \phi)$ has finitely many orbits, the possible branch lengths are $(1, \ell_1, \ell_2)$ provided that the minimum width of the branch with length 1 is 1. This is Case 2(d). Otherwise, the branch lengths must be $(1, 1, \ell)$, which is Case 2(a) or $(1, 2, \ell)$. Furthermore, by Steps 4 and 5, if the minimum width of the branch with length 1 is bigger than 2 and the minimum width of the branch with length 2 is bigger than 1, we must have $\ell \leq 4$. We conclude that  either the branch lengths are $(1, 2, \ell)$ with $\ell \leq 4$ (which is Case 2(b)); or the branch lengths are $(1,2, \ell)$ with $\ell \geq 5$ provided that either the minimum width of the branch with length 1 is at most 2 or the minimum width of the branch with length 2 is at most 1 (which is case 2(c)).

This completes the argument that any tree variety with finitely many $\PGL(n)$ orbits must be among those listed in Theorem \ref{thm-finiteorbits}. Conversely, we need to show that the tree varieties listed have finitely many orbits. In fact, using the results of Magyar, Weyman and Zelevinsky, one can enumerate all the orbits.  We will now sketch the argument.

\subsection*{Step 6: The other cases have finitely many orbits}  We have already enumerated the orbits when the tree has at most 2 branches in Step 2. We may therefore assume that the three has 3 branches, hence it looks as follows. 
 
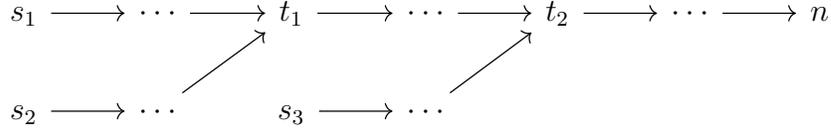
\[
\begin{tikzcd}
s_1 \arrow[r] & \cdots  \arrow[r] & t_1  \arrow[r] & \cdots \arrow[r] &  t_2 \arrow[r] & \cdots \arrow[r] & n \\
s_2 \arrow[r] &  \cdots   \arrow[ru]&  s_3 \arrow[r]& \cdots \arrow[ru]   & & 
\end{tikzcd}
\]
Two of the branches emanate from a vertex $t_1$ and the third branch emanates from a vertex $t_2$, which may coincide with $t_1$. 

If $t_1$ equals $t_2$, then let $d= d(t_1)$. The tree variety $F(T, \phi)$ maps to the partial flag variety $F(T_{\leq d}, \phi_{\leq d})$ with fiber a tree variety $F(T^1, \phi^1)$ by Proposition \ref{prop-buildinductive}. The tree $T_{\leq d}$ is a chain. Consequently, by Example \ref{ex-flag}, $F(T_{\leq d}, \phi_{\leq d})$ is a partial flag variety, hence it is homogeneous. Hence, the orbits of $\PGL(n)$ on $F(T, \phi)$ are determined by  the orbits of $\PGL(n)$ on $F(T^1, \phi^1)$. The tree variety $F(T^1, \phi^1)$ is a product of three partial flag varieties in $U_{t_1}$ by Example \ref{ex-productofflags}. In this case, the orbits are listed in \cite[Theorem 2.9]{MWZ}.

Now assume $t_1$ is different from $t_2$. Then consider the forgetful morphism forgetting all the vertices in the branches emanating from $t_1$. The image of this forgetful morphism is the tree variety corresponding to the tree
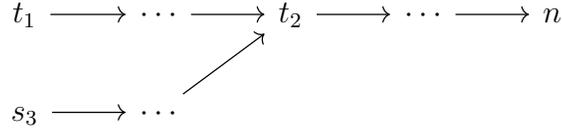
\[
\begin{tikzcd}
 t_1  \arrow[r] & \cdots \arrow[r] &  t_2 \arrow[r] & \cdots \arrow[r] & n \\
  s_3 \arrow[r]& \cdots \arrow[ru]   & & 
\end{tikzcd}
\]
By Step 2, this tree variety has finitely many orbits determined by specifying the dimension of intersections of the vector spaces in the two branches. Fix one of the orbits. The intersection of the linear spaces parameterized by the branch with leaf $s_3$  produces a partial flag in $U_{t_1}$ of specified dimension. Observe that the length of this partial flag is at most the length of the original flag (but maybe shorter if some of the linear spaces are disjoint from $U_{t_1}$). Moreover, the minimum width can only decrease. Consequently, the three partial flags in $U_{t_1}$ have finitely many orbits under the $\PGL(\phi(t_1))$-action and the orbits are listed in \cite[Theorem 2.9]{MWZ}.
\end{proof}

\section{Tree varieties with dense $\PGL(n)$ orbits}\label{sec-dense}
In this section, we discuss tree varieties that have a dense $\PGL(n)$ orbit.  Recall the following definitions from the introduction. A tree variety $F(T, \phi)$  is {\em dense} if $\PGL(n)$ has a dense orbit. Otherwise,  $F(T, \phi)$ is {\em sparse}. The tree variety $F(T, \phi)$ is trivially sparse if some vertex $v \in V(T)$ violates the inequality in Lemma \ref{lem-stilltrivial}.







We begin with an  example which will be central to our future discussion.

\begin{proposition}\label{exprop-2step}
Let $F(k_1, k_2; n)$ be a two-step flag variety. 
\begin{enumerate}
\item The triple product $F(k_1, k_2, n)^3$ is trivially sparse if and only if $n$ is divisible by 3 and $k_1 = \frac{n}{3}$ and $k_2 = \frac{2n}{3}$.
\item If  $k_1 + k_2 = n$, then   $F(k_1, k_2, n)^3$ is sparse. 
\end{enumerate}
\end{proposition}
\begin{proof}
Write $k_1 = \frac{n}{3} + x$ and $k_2 = \frac{2n}{3}+y$. Then $$n^2  - 3 \dim(F(k_1, k_2, n)) =  3x^2 - 3xy + 3y^2 = 3(x-y)^2 + 3 xy.$$ This quantity is strictly positive unless $x=y=0$. Since it is an integer, if not zero, it has to be at least 1. We conclude that the action of $\PGL(n)$ on $F(k_1, k_2, n)^3$ is trivially sparse if and only if $k_1 = \frac{n}{3}$ and $k_2 = \frac{2n}{3}$. This proves part (1).

Now assume that $k_1 + k_2 = n$. Consider the three partial flags $F_{i, k_1} \subset F_{i, k_2}$. For general choices, $F_{2, k_2}$ and $F_{3, k_2}$ intersect $F_{1, k_2}$ in two general subspaces $U_2$ and $U_3$ of dimension $2k_2 - n$. The span of $F_{2, k_1}$ and $F_{3, k_1}$ intersect $F_{1, k_2}$ in a general subspace $U_4$ of dimension $k_1$. If $F(k_1, k_2, n)^3$ is dense, then $\PGL(n)$ has to act on configurations of $F_{1, k_1}, U_2, U_3, U_4$ with a dense orbit. Since $$k_1 + \sum_{i=2}^4 \dim(U_i) = 2k_1 + 4k_2 - 2n = 2k_2,$$ by \cite[Theorem 5.1]{CHZ}, this action does not have a dense orbit. We conclude that $F(k_1, k_2, n)^3$ is sparse if $k_1 + k_2 = n$. This proves part (2).
\end{proof}

By Proposition \ref{prop-compare}, we deduce the following corollary.

\begin{corollary}\label{cor-addton}
Assume that there are two indices $i \not= j$ such that $k_i + k_j = n$. Then $F(k_1, \dots, k_r; n)^3$ is sparse.
\end{corollary}

We will next show that classifying dense tree varieties with three leaves reduces to classifying dense products of three partial flag varieties. We first introduce some notation.

\begin{notation}
Let $0< k_1 < \cdots < k_r<n$ be an increasing  sequence of positive integers less than $n$. We will denote the sequence by $k_{\bullet}$. For ease of notation, we set $k_0 =0$ and $k_{r+1} = n$.  Let $n' \leq n$ and set $d= n- n'$. Let $j$ be the index such that $k_{j-1} \leq d < k_j$. Then for $1 \leq i \leq r-j+1$, set $k_i' = k_{j+i -1} - d$.   Given a sequence $k_{\bullet}$ and an integer $n' < n$, we will call the sequence $k_{\bullet}'$  {\em the sequence derived from $k_{\bullet}$ with respect to $n$ and $n'$}.   Notice that this only depends on $n-n'$, so if we do not wish to emphasize $n$ and $n'$, we will sometimes say {\em the sequence derived from $k_{\bullet}$ with respect to $d$}. Given a vector space $W$ of dimension $n'$ and a general partial flag $F_{\bullet}$ in an $n$-dimensional vector space with dimensions $k_{\bullet}$,  the sequence $k_{\bullet}'$ denotes the dimension vector of the partial flag in $W$ obtained by $F_{\bullet} \cap W$.
\end{notation}

\begin{proposition}\label{prop-denseredprod}
Let $T$ be the following tree:
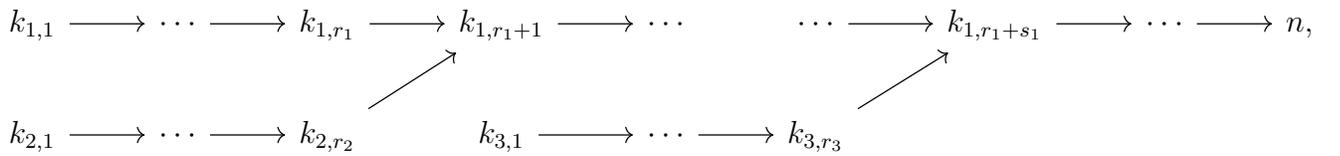
\[
\begin{tikzcd}
k_{1, 1} \arrow[r] & \cdots  \arrow[r]&  k_{1, r_1} \arrow[r] &  k_{1, r_1 + 1} \arrow[r] &  \cdots &  \cdots \arrow[r] & k_{1, r_1 + s_1} \arrow[r] & \cdots \arrow[r]&  n, \\
k_{2, 1} \arrow[r]  &  \cdots \arrow[r] &  k_{2,r_2} \arrow[ru]&  k_{3,1} \arrow[r]& \cdots \arrow[r]& k_{3, r_3} \arrow[ru]   & & 
\end{tikzcd}
\]
Set $k_{1, r_1 + 1} = m'$ and $k_{1, r_1 + s_1} =m$. Let $k_{3, \bullet}'$ be the sequence derived from $k_{3, \bullet}$ with respect to $m$ and $m'$. Then the tree variety $F(T, \phi)$ is dense if and only if 
$$\prod_{i=1}^2 F(k_{i,1}, \dots, k_{i, r_i}; m') \times F(k_{3,1}', \dots, k_{3, r_3 - j+1}', m')$$ is dense.
\end{proposition}

\begin{proof}
Let $T'$ be the tree obtained from $T$ by deleting all the vertices that have a directed path to $k_{1, r_1 + 1}$. Then $T'$ has two leaves, hence by Theorem \ref{thm-finiteorbits} it has finitely many orbits. In particular, the orbit where the two partial flags  emanating from the vertex marked $k_{1, r_1 + s_1}$ are  transverse is dense. Note that any other orbit has strictly smaller dimension and cannot contain a point of this locus in its closure by the semi-continuity of the dimension of intersections. In this orbit, the intersection of the flag $F_{3, k_{\bullet}}$ with the vector space $U_{1, k_{r_1+1}}$ has dimension vector $k_{3, \bullet}'$, the sequence derived from $k_{3, \bullet}$ with respect to $m$ and $m'$. Hence, if $F(T, \phi)$ is dense, then $\prod_{i=1}^2 F(k_{i,1}, \dots, k_{i, r_i}; m') \times F(k_{3,1}', \dots, k_{3, r_3 - j+1}', m')$ is dense.

Conversely, suppose that $\prod_{i=1}^2 F(k_{i,1}, \dots, k_{i, r_i}; m') \times F(k_{3,1}', \dots, k_{3, r_3 - j+1}', m')$ is dense. Observe that for a point in the dense orbit the flags must be pairwise transverse.  Moreover, for a transverse pair of flags, there exists a third flag so that the triple is in the dense orbit. Make the convention that $k_{2, r_2+1} = k_{1, r_1 +1}$ and $k_{r_3+1} = k_{1, r_1 + s_1}$. Fix indices so that  $k_{1, r_1 + s_1 + t_1} \to n$ is an edge in the tree.   If we omit the branch consisting of vertices $k_{2,i}$, then the resulting tree has 2 branches and has finitely many orbits. The orbit where the linear spaces are as transverse as possible is dense and has stabilizer of dimension 
$$n^2 - 1- \sum_{i=1}^{r_1+s_1 + t_1} k_{1, i} (k_{1, i+1} - k_{1, i}) - \sum_{i=1}^{r_3} k_{3, i} (k_{3,i+1} - k_{3,i}).$$
By assumption, a general choice of partial flag $U_{2,1} \subset \cdots \subset U_{2, r_2}$ in $U_{k_{1, r_1+1}}$ imposes the expected number of conditions $$\sum_{i=1}^{r_2} k_{2, i} (k_{2, i+1} - k_{2,i}).$$
Hence, the codimension of the stabilizer of such a point in $F(T, \phi)$ in $\PP GL(n)$ is the same as the dimension of $F(T, \phi)$. We conclude that $F(T, \phi)$ has a dense $\PP GL(n)$ orbit by Lemma \ref{lem-basiclem}.  
\end{proof}

Hence, for studying the density of tree varieties with three branches, it suffices to study the density of products of three partial varieties. We will concentrate on this problem for most of this section. Corollary \ref{cor-addton} and Proposition \ref{prop-denseredprod} give a large collection of sparse tree varieties.  It is also possible to give many examples of dense tree varieties.

\begin{lemma}\label{lem-easydense}
For a vertex $v \in T$, let $S(v)$ denote the set of vertices of $T$ which are sources of edges in $T$ with target $v$. If $$\sum_{s_i \in S(v)} \phi(s_i) \leq \phi(v)$$ for every vertex $v \in T$, then $F(T, \phi)$ is dense. In particular, $\prod_{i=1}^N F(k_{i,1}, \dots, k_{i, m_i}; n)$ is dense if $\sum_{i=1}^N k_{i, m_i} \leq n$.  
\end{lemma}
\begin{proof}
For each vertex $v$ with $\sum_{s_i \in S(v)} \phi(s_i) < \phi(v)$, form a new labeled tree $(T', \phi')$ by adding a new vertex $v'$ labeled $\phi'(v') = \phi(v) - \sum_{s_i \in S(v)} \phi(s_i)$ and a new edge $(v', v)$. By Proposition \ref{prop-compare}, if $F(T', \phi')$ is dense, so is $F(T, \phi)$. We may therefore assume that $T$ satisfies  $\sum_{s_i \in S(v)} \phi(s_i) = \phi(v)$ for every vertex $v$. Let $\ell_1, \dots, \ell_j$ be the leaves of the tree $T$. We have $\sum_{i=1}^j \phi(\ell_i) = n$. Fix a basis $e_1, \dots, e_n$ for the vector space $V$. Let $U_{\ell_i}$ be disjoint coordinate subspaces. The stabilizer is a block diagonal matrix with block sizes $\phi(\ell_i)$ for $1 \leq i \leq j$. Hence, the dimension of the stabilizer is $\sum_{i=1}^j \phi(\ell_i)^2 -1$. Since
$$\sum_{(s,t) \in E(T)} \phi(s)(\phi(t) - \phi(s)) = \sum_{v \in T} \sum_{s \in S(v)} \phi(s)(\phi(t)- \phi(s)) = \sum_{v \in T} (\phi(v)^2 - \sum_{s \in S(v)} \phi(s)^2) = n^2 -  \sum_{i=1}^j \phi(\ell_i)^2,$$ we conclude that $\dim(F(T, \phi))+ \dim (\Stab) = n^2-1$ and $F(T, \phi)$ is dense by Lemma \ref{lem-basiclem}.
\end{proof}

The following lemma gives a useful dimension reduction.

\begin{lemma}\label{lem-reduce}
Suppose $k_{i, \bullet}$ are $m$ flag vectors. Assume that $n'= \sum_{i=1}^{m-1} k_{i, r_i} \leq n < \sum_{i=1}^m k_{i, r_i}$. Let $k_{m, \bullet}'$ be the sequence derived from $k_{m,\bullet}$ with respect to $n$ and $n'$. Then 
$\prod_{i=1}^m F(k_{i, \bullet}; n)$ is dense if and only if $\prod_{i=1}^{m-1} F(k_{i, \bullet}; n') \times F(k_{m, \bullet}'; n')$ is dense.
\end{lemma}

\begin{proof}
Let $U_{i, \bullet}$ be general flags with dimension vectors $k_{i, \bullet}$ in $V$. For $1 \leq i \leq m-1$, the vector spaces $U_{i, k_{r_i}}$ span a vector space $W$ of dimension $n'$. The intersection $T_{m , \bullet}$ of $U_{m, \bullet}$ with $W$ has dimension vector $k_{m, \bullet}'$ and are general linear spaces. Hence, if $\prod_{i=1}^m F(k_{i, \bullet}, n)$ is dense, then $\prod_{i=1}^{m-1} F(k_{i, \bullet}; n') \times F(k_{m, \bullet}'; n')$ is dense.

The stabilizer of the flags $U_{i, \bullet}$ for $1 \leq i \leq m$, stabilizes $U_{i, \bullet}$ for $1 \leq i \leq m-1$, $W$ and $T_{m, \bullet}$. Hence, we get a map 
$$f: \Stab(\{U_{i, \bullet}\}_{i=1}^m, V) \to \Stab(\{U_{i, \bullet}\}_{i=1}^{m-1}, T_{m, \bullet}, W).$$
Pick a basis for $U_{m, \bullet}$ so that $U_{m, k_{i_m}}$ has basis $e_j$ for $n-k_{i_m} + 1 \leq j \leq n$. Let $e_1, \dots, e_{n'}$ be a basis for $W$. Then the matrices that map to the identity in $\Stab(\{U_{i, \bullet}\}_{i=1}^{m-1}, T_{m, \bullet}, W)$ have the form
$$\begin{pmatrix}
I_{n-k_{r_m}} & 0 & 0 \\
0 & I_{n'-n+k_{r_m}} & A \\
0 & 0 & B
\end{pmatrix}
$$
where $I_j$ denotes the $j \times j$ diagonal matrix and $A,B$ are obtained from the truncation of the first $n'+ k_{r_m} - n$ columns of a block lower triangular matrix of sizes $k_{m, 1}, \dots, k_{m, r_m}$. We conclude that 
\begin{equation}\label{eq-kernel}
\dim(\Ker(f)) \leq \sum_{i=1}^{j-1} k_{m,i} (k_{m,i} - k_{m, i-1}) + k_{m,j}(n-n'-k_{m, j-1}).
\end{equation}
Since $\prod_{i=1}^{m-1} F(k_{i, \bullet}; n') \times F(k_{m, \bullet}'; n')$ is dense, we have that 
\begin{multline}\label{eq-intheintersection}
\dim(\Stab(\{U_{i, \bullet}\}_{i=1}^{m-1}, T_{m, \bullet}, W)) = n'^2 - 1 - \sum_{i=1}^{m-1}\left( \left( \sum_{l=1}^{r_i-1} k_{i,l} (k_{i,l+1} - k_{i,l}) \right)+ k_{i,r_i}(n' - k_{i,r_i})\right) \\ - \sum_{l = j}^{r_m-1} (k_{m, l} - n + n') (k_{m, l+1} - k_{m, l}) - (k_{m, r_m} - n + n') (n - k_{m, r_m}) \\
= n^2-1 - \sum_{i=1}^m\sum_{l=1}^{r_i-1} k_{i,l}(k_{i, l+1} - k_{i,l}) - \sum_{i=1}^m k_{i,r_i}(n-k_{i,r_i}) + \sum_{l=1}^{j-1} k_{m,l}(k_{m,l+1}-k_{m,l}) -(n-n') k_{m,j}
\end{multline}
By the theorem on the dimension of fibers, we have 
$$\dim(\Stab(\{U_{i, \bullet}\}_{i=1}^m, V)) \leq \dim(\Stab(\{U_{i, \bullet}\}_{i=1}^{m-1}, T_{m, \bullet}, W)) + \dim(\Ker(f)).$$ Combining this with Equations \eqref{eq-kernel} and \eqref{eq-intheintersection} and some arithmetic, we see that 
$$\dim(\Stab(\{U_{i, \bullet}\}_{i=1}^m, V)) \leq  n^2-1 - \sum_{i=1}^m \sum_{l =1}^{r_i-1} k_{i,l} (k_{i, l+1} - k_{i,l}) - \sum_{i=1}^m k_{i, r_i}(n-k_{i,r_i}).$$
Hence, $\prod_{i=1}^m F(k_{i, \bullet}, n)$ is dense by Lemma \ref{lem-basiclem}.
\end{proof}

\begin{lemma}\label{lem-n2kreduction}
Let $n=2k_r$. Then $F(k_1, \dots, k_r; n)^3$ is dense if and only if $F(k_1, \dots, k_{r-1}; k_r)^3$ is dense.
\end{lemma}

\begin{proof}
Fix three general partial flags $W_{\bullet}^i \in F(k_1, \dots, k_r; n)$ for $1 \leq i \leq 3$. Let $Y_j'$ be the span of $W_j^2$ and $W_{k_r}^3$. Let  $Z_j'$ be the span of $W_j^3$ with $W_{k_r}^2$. Set $Y_j = Y_j' \cap W_{k_r}^1$ and $Z_j = Z_j' \cap W_{k_r}^1$. We have that $\dim (Y_j) = \dim(Z_j) = j$. Observe that a general pair of partial flags $(Y_{\bullet}, Z_{\bullet}) \in F(k_1, \dots, k_{r-1}; k_r)^2$ occurs this way. To see this, fix two general $k_r$-dimensional subspaces $W_{k_r}^2$ and $W_{k_r}^3$. We can recover the partial  flag $W_{\bullet}^2$ by taking the span of $Y_{\bullet}$ with $W_{k_r}^3$ and intersecting with $W_{k_r}^2$. We can recover the partial  flag $W_{\bullet}^3$ by taking the span of $Z_{\bullet}$ with $W_{k_r}^2$ and intersecting with $W_{k_r}^3$. Now the construction yields back the partial flags $Y_{\bullet}$ and $Z_{\bullet}$. Hence, if $F(k_1, \dots, k_r; n)^3$ is dense, then $F(k_1, \dots, k_{r-1}; k_r)^3$ is dense.

Conversely, suppose $F(k_1, \dots, k_{r-1}; k_r)^3$ is dense. We get a homomorphism 
$$f: \Stab(W_{\bullet}^1, W_{\bullet}^2, W_{\bullet}^3, \CC^n) \to \Stab(W_{\bullet}^1, Y_{\bullet}, Z_{\bullet}, W_{k_r}^1).$$
Choose a basis for $\CC^n$. Set $W_{k_r}^1$ to be the span of $e_i$ with $1 \leq i \leq k_r$. Set  $W_{k_r}^1$ to be the span of $e_i$ with $k_r+ 1 \leq i \leq n$. Finally, set $W_{k_r}^3$ to be the span of $e_i + e_{i + k_r}$ for $1 \leq i \leq k_r$. Then the stabilizer of the three subspaces has the form $$\begin{pmatrix}
A & 0 \\
0 & A 
\end{pmatrix},
$$
where $A$ is a $k_r \times k_r$ invertible matrix. Hence, the kernel of the map $f$ is trivial. We conclude that 
$$\dim(\Stab(W_{\bullet}^1, W_{\bullet}^2, W_{\bullet}^3, \CC^n)) \leq \dim(\Stab(W_{\bullet}^1, Y_{\bullet}, Z_{\bullet}, W_{k_r}^1)) = k_r^2 -1 - 3 \sum_{i=1}^{r-1} k_i (k_{i+1} - k_i).$$
On the other hand, 
$$\dim(\Stab(W_{\bullet}^1, W_{\bullet}^2, W_{\bullet}^3, \CC^n)) \geq n^2 -1 - 3 \sum_{i=1}^{r} k_i (k_{i+1} - k_i).$$
Since $n=2k_r$,
$$n^2 -1 - 3 \sum_{i=1}^{r} k_i (k_{i+1} - k_i)= k_r^2 -1 - 3 \sum_{i=1}^{r-1} k_i (k_{i+1} - k_i)$$ and we have equality every where. We conclude that $F(k_1, \dots, k_r; n)^3$ is dense.
\end{proof}

\begin{proposition}\label{prop-double}
Assume that $2 k_r \leq n$ and $2 k_i \leq k_{i+1}$ for   $2 \leq i  \leq r-1$.  Then $F(k_1, \dots, k_r; n)^3$ is dense.
\end{proposition}

\begin{proof}
We will prove the proposition by induction on $r$.  If $r=1$, the proposition is true by \cite[Theorem 5.1]{CHZ}. Now suppose the proposition is true up to $r-1$. If $3k_r \leq n$, then  Lemma \ref{lem-easydense} implies that  $F(k_1, \dots, k_r; n)^3$ is dense. We may therefore assume that $2k_r \leq n < 3k_r$. 

If $2k_r = n-m$ with $m >0$, then let $k_{\bullet}'$ be the sequence derived from $k_{\bullet}$ with respect to $m$. By applying Lemma \ref{lem-reduce} three times, the density of $F(k_{\bullet}, n)^3$ is equivalent to the density of $F(k_{\bullet}'; n-3m)^3$. Observe that $2(k_r-m) = n-3m$ and $2(k_i - m) = 2k_i - 2m < k_{i+1} - m$. Hence, $k_{\bullet}'$ still satisfies the assumptions of the proposition.  We therefore reduce to the case when $2k_r = n$. 

By Lemma \ref{lem-n2kreduction},  $F(k_1, \dots, k_r; 2k_r)^3$ is dense  if and only if $F(k_1, \dots, k_{r-1}; k_r)^3$ is dense. The latter satisfies the assumptions of the proposition and has one fewer steps. Hence, by induction, it is dense. This concludes the proof of the proposition.
\end{proof}

For our inductive arguments, we will need a technical lemma. Let $V$ be an $n$-dimensional vector space. For $1 \leq i \leq 3$, let $U_i \subset T_i$ be three two-step flags of dimensions $u_i < t_i$  in $V$. Let $W \subset V$ be a subspace of dimension $w$ containing $U_2$ and $T_3$. Let $X$ denote the variety which parameterizes such configuration of subspaces of $V$. Let $m = n-w$ and assume that $m \leq u_1$ and $m \leq t_2 - u_2$.
Let $U_1' \subset T_1'$ be a two-step flag in $W$ of dimension $u_1 - m, t_1 - m$. Let $U_2' \subset T_2'$ be a two-step in $W$ of dimension $u_2, t_2 - m$ and let $U_3' \subset T_3'$ be a two-step flag in $W$ of dimension $u_3, t_3$. Finally, let $W'$ be a linear subspace of dimension $u_1 + t_2 - m$ containing $U_1'$ and $T_2'$. Let $Y$ be the variety which parameterizes such configurations of subspaces of $W$.

\begin{lemma}\label{lem-techred}
With this notation, assume that  $w \geq u_2 + t_3, \ v_1 + v_2 > n, \ u_1 + v_2 < n.$ Then:
\begin{enumerate}
\item The variety $X$ is irreducible of dimension 
$$ \sum_{i=1}^3(u_i (t_i-u_i) + t_i(n-t_i)) + (w-u_2-t_3)m.$$
\item The variety $Y$ is irreducible of dimension 
$$\sum_{i=1}^3(u_i (t_i-u_i) + t_i(n-t_i)) + (w-u_2-t_3)m - 2mw - m^2.$$
\item The $\PGL(n)$ action on $X$ has a dense orbit if and only if the $\PGL(w)$ action on $Y$ has a dense orbit.
\end{enumerate}
\end{lemma}

\begin{proof}
We first observe that the varieties $X$ and $Y$ parameterizing the specified configurations are irreducible varieties. If we omit $T_2$, then we obtain the tree variety associated to the tree
$$
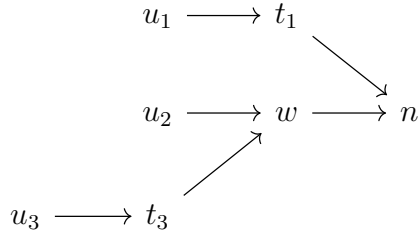
\begin{tikzcd}
 & u_1 \arrow[r] & t_1  \arrow[rd] &  \\
 &  u_2 \arrow[r] & w \arrow[r] & n \\
u_3 \arrow[r] & t_3 \arrow[ru] & 
\end{tikzcd} $$
The choice of $T_2$ containing $U_2$ realizes $X$ as a Grassmannian bundle $G(t_2-u_2, V/U_2)$ over this tree variety. By Theorem \ref{thm-irredanddim}, $X$ is irreducible and its dimension is as claimed. A similar argument shows that $Y$  irreducible of the claimed dimension. Omitting $T_3'$ gives rise to the tree variety associated to the tree
$$
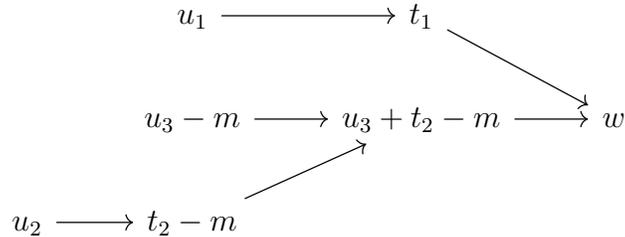
\begin{tikzcd}
 & u_1 \arrow[r] & t_1  \arrow[rd] &  \\
 &  u_3-m \arrow[r] & u_3 + t_2 - m \arrow[r] & w \\
u_2 \arrow[r] & t_2-m \arrow[ru] & 
\end{tikzcd} $$
The choice of $T_3'$ containing $U_3'$ realizes $Y$ as a Grassmannian bundle $G(t_3-u_3, W/U_2')$  over this tree variety. By Theorem \ref{thm-irredanddim}, $Y$ is irreducible of the claimed dimension.

There is a rational map $X \dashrightarrow Y$ given by setting $$U_i' = U_i \cap W, \quad T_i' = T_i \cap W, \quad W' = \overline{U_1T_2} \cap W.$$ A general configuration in $W$ occurs as the intersection of a configuration in $V$ with $W$. Hence, if the configuration in $V$ has a dense orbit, then the configuration in $W$ has a dense orbit. We need to prove the converse. Given a general configuration, we obtain a homomorphism 
$$f: \Stab(U_i, T_i, W; V) \to \Stab(U_i', T_i', W'; W).$$ First observe that the kernel of $f$ is trivial. To see this, we may choose a basis $e_i$, $1 \leq i \leq n$, for $V$ so that $W$ is spanned by $e_i$ with $1 \leq i \leq w$ and $T_2$ is spanned by $e_i$ with $n-t_2 + 1 \leq i \leq n$. Finally, we may choose $T_1$ to be spanned by $e_i + e_{w+i}$ for $1 \leq i \leq n $

Hence, $$\dim (\Stab(U_i, T_i, W; V)) \leq \dim (\Stab(U_i', T_i', W'; W))$$ Since $Y$ has a dense orbit,  the dimension of $\Stab(U_i', T_i', W'; W)$ is $$w^2 -1+ 2mw + m^2 - \sum_{i=1}^3(u_i (t_i-u_i) + t_i(n-t_i)).$$
Since $n=w+m$, we have  $$n^2 -1  - \sum_{i=1}^3 u_i (t_i - u_i) - \sum_{i=1}^3 t_i(n-t_i) = w^2 -1 - \sum_{i=1}^3 u_i (t_i - u_i) - \sum_{i=1}^3 t_i(n-t_i) + 2mw + m^2.$$ Since the latter is the dimension of the stabilizer of a general point in $Y$ and bounds the dimension of the stabilizer of a general point in $X$, we conclude that if $Y$ has a dense orbit, then so does $X$.
\end{proof}

\begin{theorem}\label{thm-2stepmain}
Let $F(k_1, k_2; n)$ be a two-step partial flag variety. Then $F(k_1, k_2; n)^3$ is sparse if and only if $k_1 + k_2 = n$.
\end{theorem}
\begin{proof}
By Proposition \ref{exprop-2step}, we know that if $k_1 + k_2 = n$, then $F(k_1, k_2; n)^3$ is sparse. We need to show that if $k_1 + k_2 \not= n$, then $F(k_1, k_2; n)^3$ is dense.

By replacing $F(k_1, k_2;n)$ with $F(n-k_2, n-k_1; n)$ if necessary, we may assume that $k_1 + k_2 < n$.
If $3k_2 \leq n$, then $F(k_1, k_2; n)^3$ is dense  by Lemma \ref{lem-easydense}. 

If $2k_2 \leq n < 3k_2$, then let $m = n -2k_2$. If $m \geq k_1$, then we apply Lemma \ref{lem-reduce} three times. First, the density of $F(k_1, k_2; n)^3$ is equivalent to the density of $F(k_1, k_2; 2k_2)^2 \times G(k_2-m; 2k_2)$. The density of the latter is in turn equivalent to the density of $F(k_1, k_2; 2k_2-m) \times G(k_2-m; 2k_2 -m)^2$, which is equivalent to the density of $G(k_2-m; 2k_2-2m)^3$. Since the latter is dense by \cite[Theorem 5.1]{CHZ}, we conclude that $F(k_1, k_2; n)^3$ is dense. If $m < k_1$, then  by applying Lemma \ref{lem-reduce} three times, the density of $F(k_1, k_2; n)^3$ is equivalent to the density of $F(k_1-m, k_2-m; 2k_2-2m)^3$. By Lemma \ref{lem-n2kreduction}, the density of the latter is equivalent to the density of $G(k_1-m; k_2-m)^3$. Since the latter is dense, we conclude that $F(k_1, k_2; n)^3$ is dense in this case.

Finally, we may assume that $k_2 < n < 2k_2$.  If $2k_1 + k_2 < n$, then let $m=n-2k_1-k_2$. By applying Lemma \ref{lem-reduce} three times, the density of $F(k_1, k_2; n)^3$ is equivalent to the density of $F(k_1, k_2-2m; n-3m)^3$. Hence, we reduce to the case $2k_1 + k_2 \geq n$.

If $2k_1 + k_2 \geq n$, let $U_i \subset T_i$ for $1 \leq i \leq 3$ be  three general partial flags of type $k_1, k_2$. Let $m = n - k_1 - k_2$. Let $W_{i,j}$ denote the span of $U_i$ and $T_j$. We will apply Lemma \ref{lem-techred} three times to reduce the density of $F(k_1, k_2; n)^3$ to that of $F(k_1-m, k_2-2m; n-3m)^3$.   First apply Lemma \ref{lem-techred}, setting $W= W_{2,3}$. Denote the intersection of a linear space with $W$ with a prime.  Then $F(k_1, k_2; n)^3$ is dense if and only if the configuration $(U_i', T_i', W_{1,2}')$ is dense in $W$. Now apply Lemma \ref{lem-techred} setting $W = W_{1,2}'$. Denote the intersections of the vector spaces with $W_{1,2}'$ with double primes. Set $W_{3,1}'' : = \overline{U_3'T_1'} \cap W_{1,2}'$. Then the configuration $(U_i', T_i', W_{1,2}')$  in $W$ is dense if and only if the configuration $(U_i'', T_i'', W_{3,1}'')$ is dense in $W_{1,2}'$. Finally, we apply Lemma \ref{lem-techred} by setting $W=W_{3,1}''$. Denote the intersections of the vector spaces with $W_{3,1}''$ by triple primes. Notice that $U_2$ and $T_3''$ span $W_{3,1}''$. We conclude that the configuration $(U_i'', T_i'', W_{3,1}'')$ is dense in $W_{1,2}'$ if and only if the configuration $(U_i''', T_i''')$ is dense in $W_{3,1}''$. We have thus reduced the density of $F(k_1, k_2; n)^3$ to that of $F(k_1-m, k_2-2m; n-3m)^3$. Notice that $k_1 - m + k_2 - 2m < n -3m$ by assumption. If $2k_2 - 4m \leq n-3m$, we are done by the previous cases. Otherwise, we can continue reducing the size of $k_1$ and $k_2$ by $m$ and $2m$, respectively. Since this cannot continue indefinitely, we conclude that  $F(k_1, k_2; n)^3$ is dense. This concludes the proof of the theorem.

\end{proof}

We conclude with a few remarks about the action of $\PGL(n)$ on products of Grassmannians. Classifying the products of at least 5 Grassmannians with dense orbit is a hard problem. However, one can say a little more about certain families of such products.

\begin{lemma}
The action of $\PGL(n)$ on $\prod_{i=1}^m G(k_i, n)$ has a dense orbit if $$\sum_{i=1}^{m-1} k_i \leq n \quad \mbox{and} \quad k_m \leq n - \sum_{i=1}^{m-1} k_i + \min_{1 \leq i \leq m-1} k_i.$$
\end{lemma}
\begin{proof}
For simplicity set $s= n - \sum_{i=1}^{m-1} k_i$.
Fix a basis $e_j$, $1 \leq j \leq n$, of $V$. For $1 \leq i \leq m-1$, let $W_i$ be the vector space spanned by $e_j$ with $1+ \sum_{l=1}^{i-1} k_l \leq j \leq  \sum_{l=1}^{i} k_l$. Let $W_m$ be the vector space spanned by $e_j$ for $1+ \sum_{l=1}^{m-1} k_l \leq j \leq n$ and $e_j + e_{j+k_1} + e_{j + k_1 + k_2} + \cdots + e_{j + \sum_{i=1}^{m-1} k_i}$ for $1 \leq j \leq k_m -s$. Then the stabilizer of these vector spaces have the form
$$\begin{bmatrix}
A & B_1 & 0& 0 & \cdots & D \\
0 & C_1 & 0& 0  & \cdots  & 0 \\
0 & 0 & A& B_2 & \cdots & D \\
0 & 0 & 0 & C_2 & \cdots & 0 \\
& \cdots & & & & \\
0 & 0 & 0 & 0 &  \cdots & E
\end{bmatrix}$$
This stabilizer of this configuration has dimension $$(k_m-s)^2 + \sum_{i=1}^{m-1} k_i (k_i -k_m + s) + s^2 + s(k_m -s) - 1.$$
Since $n= \sum_{i=1}^{m-1} k_i + s$, we conclude that this quantity is equal to 
$$\sum_{i=1}^m k_i ^2 + n (s- k_m ) -1.$$ Observe that this quantity is also equal to 
$$n^2 -1 - \sum_{i=1}^m k_i (n-k_i) = n^2 -1 + \sum_{i=1}^m k_i^2 - n(n-s+k_m) =  \sum_{i=1}^m k_i^2 + n(s-k_m) -1.$$ Since $\dim(\Stab(W_i; V))= \dim(\PGL(n)) - \dim (\prod_{i=1}^m G(k_i, n)),$ we conclude that $\PGL(n)$ has a dense orbit on $\prod_{i=1}^m G(k_i, n)$.
\end{proof}

\begin{example}
The assumption 
$k_m \leq n - \sum_{i=1}^{m-1} k_i + \min_{1 \leq i \leq m-1} k_i$ cannot be weakened in general.  For example, set $n=8$ and $k_1 = k_2 = 1, k_3=k_4 = 2$ and $k_5=3$. Then $\sum_{i=1}^4 k_i = 8 = n$. However, $3>  n - \sum_{i=1}^{m-1} k_i + \min_{1 \leq i \leq m-1} k_i =1$. For $1 \leq i \leq 5$, let $W_i$ be a general linear space of dimension $k_i$. Let $I$ be any three element subset of $\{1, 2, 3, 4\}$ and let $j$ be the element in the complement of $I$. Then the span of $W_i$ for $i \in I$ intersects $W_5$ in a subspace of dimension $3-k_j$. In this way, we get 4 general subspaces of dimensions $1,1,2,2$ in $W_5$. Since $\PGL(3)$ does not have a dense orbit on this configuration, the original configuration does not have a dense orbit.
\end{example}

\begin{theorem}$\label{thm:abcde}$
Let $\underline{\textbf{d}}=(d_1,d_2,d_3,d_4,n-d_5;n)$ be a dimension vector such that $n\geq d_1+d_2+d_3+d_4$ and $d_5\geq d_4 \geq d_3 \geq d_2 \geq d_1$. Then $\underline{\textbf{d}}$ is  dense if and only if $d_1+d_2+d_3+d_4\not= 2d_5$
\end{theorem}
\begin{proof}

Let $V$ be an $n$ dimensional vector space. By Proposition $\ref{prop-complement}$ we can consider the vector spaces $V_{d_1},V_{d_2},V_{d_3},V_{d_4},W$ with the corresponding dimensions $n-d_1,n-d_2,n-d_3,n-d_4,d_5$ respectively. Now consider the following group homomorphism that is constructed by the restriction map:
$$f:\Stab(V_{d_1},V_{d_2},V_{d_3},V_{d_4},W;V) \to \Stab(V_{d_1}\cap W,V_{d_2} \cap W,V_{d_3} \cap W,V_{d_4}\cap W;W)  $$
Hence we have 
$$\dim \Stab(V_{d_1},V_{d_2},V_{d_3},V_{d_4},W;V) = \dim \Stab(V_{d_1}\cap W,V_{d_2} \cap W,V_{d_3} \cap W,V_{d_4}\cap W;W) + \dim \ker f $$

Now we show that kernel of $f$ is trivial. For this purpose we consider an element in the preimage of the identity in $\text{Stab}(V_{d_1}\cap W,V_{d_2} \cap W,V_{d_3} \cap W,V_{d_4}\cap W;W)$. Choose coordinates such that $V_{d_1}$ be given by $x_1=x_2=..=x_{d_1} = 0$, $V_{d_2}$ be given by $x_{d_1+1}=..=x_{d_2}=0$,  $V_{d_3}$ be given by $x_{d_2+1}=..=x_{d_3}=0$ and  $V_{d_4}$ be given by $x_{d_3+1}=..=x_{d_4}=0$.
Then $\Stab(V_{d_1},V_{d_2},V_{d_3},V_{d_4};V)$ is given by the $n\times n$ block matrix:
$$\begin{bmatrix}\label{matrix:blockMat}
A & 0 & 0& 0\\
0 & B & 0& 0\\
0 & 0 & C& 0\\
0 & 0 & 0 & D
\end{bmatrix}$$
where $A,B,C,D$ are $d_1 \times d_1, d_2 \times d_2, d_3 \times d_3, d_4 \times d_4 $ matrices, respectively.
Finally we show that if this matrix acts as an identity on $W$ then it is in fact the identity. Now let $W$ be spanned by the vectors of the $\ref{matrix:blockMat}$ such that it has a basis $\{ e_1+e_{d_1+1}+e_{d_1+d_2+1}+e_{d_1+d_2+d_3+1},e_{d_1+d_2+d_3+d_4+1},e_{d_1+2}+..+e_{d_1+d_2+d_3+d_4+2},..\}$ where an element from each block of A,B,C,D is taken for each basis element and summed through all blocks and if one element does not exist, we omit that element. As we have $d_5\geq d_4$ we use all the elements in the block matrix and hence if we fix those elements and their sums with a element in $\text{Stab}(V_{d_1}\cap W,V_{d_2} \cap W,V_{d_3} \cap W,V_{d_4}\cap W;W)$, then the matrix has to be identity matrix. Hence $\Stab(V_{d_1},V_{d_2},V_{d_3},V_{d_4},W;V)$ acts with a dense orbit if and only if $\Stab(V_{d_1}\cap W,V_{d_2} \cap W,V_{d_3} \cap W,V_{d_4}\cap W;W)$ acts with a dense orbit.
\\
 Now from \cite[Theorem 5.1]{CHZ} we have $\Stab(V_{d_1}\cap W,V_{d_2} \cap W,V_{d_3} \cap W,V_{d_4}\cap W;W)$ does not act with a dense orbit if and only if $d_1+d_2+d_3+d_4=2d_5$ and the theorem follows.
\end{proof}

\bibliographystyle{plain}

\end{document}